\documentclass[a4paper]{siamart1116}

\usepackage{dsfont}
\usepackage{amsmath}
\usepackage{amssymb}
\usepackage{amsxtra}

\usepackage{amsfonts}
\usepackage{graphicx}
\usepackage{epstopdf}
\usepackage{algorithmic}
\Crefname{ALC@unique}{Line}{Lines} % <- Preamble
\usepackage{caption}
\usepackage{upquote}
\usepackage{longtable}

\ifpdf
  \DeclareGraphicsExtensions{.eps,.pdf,.png,.jpg}
\else
  \DeclareGraphicsExtensions{.eps}
\fi

\usepackage{geometry}
\numberwithin{theorem}{section}

\newcommand{\TheTitle}{On computing the Kronecker structure of polynomial and rational matrices using \textsf{J\lowercase{ulia}}}
\newcommand{\TheAuthors}{A. Varga}

\headers{\TheTitle}{\TheAuthors}

\title{{\TheTitle}}%\thanks{Submitted to the editors DATE.
\author{
  Andreas Varga\thanks{Gilching, Germany
    (\email{varga.andreas@gmail.com}).}
}
\headers{Computing the Kronecker Structure of Polynomial and Rational Matrices}{A. Varga}

\ifpdf
\hypersetup{
  pdftitle={\TheTitle},
  pdfauthor={\TheAuthors}
}
\fi

\newcommand{\be}{\begin{equation}}
\newcommand{\ee}{\end{equation}}
\newcommand{\ba}{\left [ \begin{array}}
\newcommand{\ea}{\end{array} \right ]}
\newcommand{\bea}{\begin{eqnarray}}
\newcommand{\eea}{\end{eqnarray}}

\newcommand{\finr}{{\ \hfill $\Box$}}

\newcommand{\rev}{\mathop{\mathrm{rev}}}

\newcommand{\rank}{\mathop{\mathrm{rank}}}
\newcommand{\diag}{\mathop{\mathrm{diag}}}

\begin{document}

\maketitle

% REQUIRED
\begin{abstract}
In this paper we discuss the mathematical background and the computational aspects which underly the implementation of a collection of \textsf{Julia} functions in the \textsf{MatrixPencils} package for the determination of structural properties of polynomial and rational matrices. We primarily focus on the
computation of the finite and infinite spectral structures (e.g., eigenvalues, zeros, poles) as well as the left and right singular structures (e.g., Kronecker indices), which play a fundamental role in the structure of the solution of many problems involving polynomial and rational matrices. The basic analysis tool is the determination of the Kronecker structure of linear matrix pencils using numerically reliable algorithms, which is used in conjunction with several linearization techniques of polynomial and rational matrices.  Examples of polynomial and rational matrices, which exhibit all relevant structural features, are considered to illustrate the main mathematical concepts and the capabilities of implemented tools.
\end{abstract}

% REQUIRED
\begin{keywords}
  Polynomial matrices, rational matrices, matrix pencils, descriptor systems, computational methods.
\end{keywords}

% REQUIRED
\begin{AMS}
   	26C10, 30C10, 93B10, 93B60, 93C05
\end{AMS}

\section{Introduction}

Structural properties such as eigenvalues, zeros, poles, and minimal indices, play a fundamental role in the structure of the solution of many problems involving polynomial and rational matrices. An important application domain is the solution of polynomial eigenvalue problems, where polynomial matrices arise either directly from the mathematical modelling of continuous system dynamics or represent approximations of more general nonlinear mappings leading to nonlinear eigenvalue problems (see \cite{Mack15} for a fairly complete survey of this subject). Another field of application is control system theory, where polynomial and rational matrix models play a fundamental role in the structural analysis of linear systems \cite{Rose70}.

One of the computational approaches to solve polynomial eigenvalue problems is via linearizations, where polynomial matrices of arbitrary degree are replaced by first degree polynomial matrices (also called \emph{matrix pencils}) which allow to retrieve the structural feature of the original problems. The main appeal of this approach is to allow the use of well established computational techniques for matrix pencil manipulations (e.g., reduction to various Kronecker-like forms in conjunction with the QZ algorithm) to determine the involved structural elements. The most commonly used linearizations are the Frobenius companion forms \cite{DeTe14}, which can be directly built from the underlying problem data.

For rational matrices (in particular also for polynomial matrices), alternative linearizations are structured matrix pencils, also called system matrix pencils  \cite{Rose70}, which share the same pole-zero and singular structures with the original rational matrices. Using these linearizations for pole-zero computations or for the determination of the singular structures of rational matrices requires additionally the determination of linearization with special features (e.g., of least dimension). Therefore, minimal realization methods, as those employed in control system theory, are basic computational ingredients to determine least order linearizations.

In this article we present the basic concepts to characterize the structural properties of polynomial and rational matrices such as finite and infinite eigenvalues, minimal indices, zeros and poles, and  discuss these concepts also in the particular case of first degree polynomial matrices (i.e., matrix pencils). Numerically reliable matrix pencil reduction techniques play a central role in the determination of these properties and therefore they form the basic numerical ingredients for the investigation of structural features of polynomial matrices via suitable linearizations. Three classes of linearizations are discussed for a given polynomial matrix (companion forms, pencil based system matrix, and descriptor system matrix) and the correspondences between the properties of the original polynomial matrix and its linearizations are described. Additionally, we describe a general linearization technique of structured polynomial system matrices of arbitrary degree. For rational matrices, pencil based and descriptor system based linearizations form the basis of computational methods for structural analysis. We present succinctly the newly implemented collection of software tools for the \textsf{Release v1.1} of the \textsf{Julia} package \textsf{MatrixPencils}, which cover the  computation of structural elements of polynomial and rational matrices and related computations as described in this paper.
The main mathematical concepts and the capabilities of implemented tools  are illustrated using an example of a polynomial matrix employed in \cite{Door83a} and a rational matrix derived from this example,  which exhibit all relevant structural features.

\section{Polynomial matrices}
Let $\mathds{F}(\lambda)$ be the set of rational functions with coefficients in the field $\mathds{F}$ with indeterminate $\lambda$, and let $\mathds{F}[\lambda]$ be the set of polynomials with coefficients in the field $\mathds{F}$. We denote $\overline{\mathds{F}}$ the algebraic closure of $\mathds{F}$.  The most usual cases are when either $\mathds{F} = \mathds{R}$, the set of real numbers, or  $\mathds{F} = \mathds{C}$, the set of complex numbers. Since polynomials can be assimilated with special rational functions with 1 as denominator, $\mathds{F}[\lambda] \subset \mathds{F}(\lambda)$. It is easy to show that $\mathds{F}(\lambda)$ is closed under the addition and multiplication operations.  Both operations are associative and commutative, the multiplication is distributive over addition,
and each operation possesses an identity element in $\mathds{F}(\lambda)$. Finally, there exist inverses for
all elements under addition and for all nonzero elements under
multiplication. Therefore, the set $\mathds{F}(\lambda)$ forms a \emph{field}. The subset of polynomials $\mathds{F}[\lambda]$ forms only a \emph{ring} (more exactly, an Euclidean domain with identity), because the only invertible elements in $\mathds{F}[\lambda]$ are the nonzero elements of $\mathds{F}$, which are thus the \emph{units} of the ring.

Let $P(\lambda) \in \mathds{F}[\lambda]^{m\times n} $ be a $m\times n$ \emph{polynomial matrix}  defined as
\be\label{polmat} P(\lambda) = \ba{ccc} p_{11}(\lambda) & \cdots & p_{1n}(\lambda) \\
\vdots & \ddots & \vdots \\
p_{m1}(\lambda) & \cdots & p_{mn}(\lambda) \ea \, , \ee
where each $p_{ij}(\lambda)$ is a polynomial of the form
\be\label{polij} p_{ij}(\lambda) = a_k\lambda^k+a_{k-1}\lambda^{k-1}+\cdots + a_1\lambda+a_0  \ee
with coefficients in $\mathds{F}$. Polynomial row vectors, column vectors and even scalar polynomials can be associated to particular polynomial matrices with $m = 1$, $n = 1$ or $m = n = 1$, respectively.
The degree $d$ of $P(\lambda)$ is the largest degree of the polynomial entries of $P(\lambda)$
\[ d = \deg P(\lambda) := \max_{i,j} \deg p_{ij}(\lambda) . \]
Polynomials as in (\ref{polij}), with $a_k = 1$, are called \emph{monic polynomials}.

If $k \geq d$, $P(\lambda)$ can alternatively be written as a \emph{grade} $k$ matrix polynomial
\be\label{matpol} P(\lambda) = \sum_{i = 0}^k \lambda^i P_i  \ee
with  $P_i \in \mathds{F}^{m\times n}$, for $i = 0, 1, \ldots, k$. For this representation, the degree is simply the largest index $d$ for which $P_d \neq 0$. A polynomial matrix with $k = d = 1$ is called a \emph{matrix pencil}, or simply a \emph{pencil}.

\emph{Remark. }The choice of grade $k$ is a matter of convenience and depends on the intended application. For example, it allows to view a constant polynomial matrix $P(\lambda) = A$ as $A$ (i.e., a constant matrix), or $A + \lambda 0 $ (i.e., a pencil), or even as $A + \lambda 0 + \lambda^2 0 + \ldots + \lambda^k 0$ (i.e., a polynomial matrix of grade $k$). For a pertinent discussion of this matter, see \cite{DeTe14}. \finr

The following definitions are straightforward extensions of familiar notions for constant matrices.  $P(\lambda)$ is called \emph{regular} if $m = n$ and $\det P(\lambda) \not\equiv 0$. Otherwise, $P(\lambda)$ is called \emph{singular}. Equivalently, $P(\lambda)$,  viewed as a rational matrix with entries in the field $\mathds{F}(\lambda)$, is regular if $P(\lambda)$ is invertible (the inverse is however not a polynomial matrix in general). The normal rank of $P(\lambda)$, denoted $\rank P(\lambda)$,  is the size of the
largest non-identically-zero minor of $P(\lambda)$. Equivalently, the normal rank of $P(\lambda)$, viewed
as a rational matrix, is the number of linearly independent rows or columns of $P(\lambda)$. A regular polynomial matrix $P(\lambda)$ is called \emph{unimodular} if $\det P(\lambda)$ is a constant (i.e., independent of $\lambda$), or, equivalently, $P(\lambda)$ has an inverse that is also a polynomial matrix.

Two main structural properties of a polynomial matrix $P(\lambda)$ are its eigenvalue structure and its singular structure.
In what follows, we address these aspects using both linear algebra results as well well as control system theory results.

The \emph{eigenvalue structure} concerns with the eigenvalues of the polynomial matrix $P(\lambda)$, which are those values of $\lambda$ for which the equation
\[ P(\lambda) x = 0 \]
has nonzero solutions $x$. For example, if $P(\lambda)$ is regular, then the finite eigenvalues are simply the roots of $\det P(\lambda)$. If $\rank P_d = n$, this is a polynomial of degree $nd$ and all eigenvalues of $P(\lambda)$ are finite. If $\rank P_d < n$, then $\det P(\lambda)$ is a polynomial of degree say $q < nd$ and, therefore $P(\lambda)$ has $q$  finite eigenvalues and $nd-q$ infinite eigenvalues. In what follows, we give the precise definitions of eigenvalues using the Smith canonical form of polynomial matrices.

\begin{theorem}[Smith form]\label{SF}
Let $P(\lambda)$ be an $m\times n$ polynomial matrix of rank $r$ with coefficients in  $\mathds{F}$. Then, there exist unimodular polynomial matrices $U(\lambda) \in \mathds{F}[\lambda]^{m\times m}$ and $V(\lambda) \in \mathds{F}[\lambda]^{n\times n}$ such that
 \be\label{SmithForm} D(\lambda) := U(\lambda)P(\lambda)V(\lambda) =
 \ba{ccc|c} d_1(\lambda) &&& \\ & \ddots && 0_{r,n-r} \\ &&d_r(\lambda) & \\ \hline & 0_{m-r,r} && 0_{m-r,n-r} \ea , \ee
 where $d_1(\lambda)$, $\cdots$, $d_r(\lambda)$ are monic polynomials in $\mathds{F}[\lambda]$ such that  $d_i(\lambda)$ divides $d_{i+1}(\lambda)$ for $i = 1, \ldots, r-1$. Moreover, $D(\lambda)$ is unique and is called the Smith canonical form of $P(\lambda)$.
\end{theorem}
The monic polynomials $d_1(\lambda)$, $\cdots$, $d_r(\lambda)$ are called the \emph{invariant polynomials} of $P(\lambda)$.
The (finite) eigenvalues of $P(\lambda)$ are the totality of (finite) zeros (roots) of all invariant polynomials.

For each distinct eigenvalue $\lambda_0 \in \overline{\mathds{F}}$, we can express each $d_i(\lambda)$ in a factored form as $d_i(\lambda) = (\lambda-\lambda_0)^{\alpha_i}p_i(\lambda)$ with $p_i(\lambda_0) \neq 0$, where $\alpha_i \geq 0$ is called the $i$-th \emph{partial multiplicity} of $\lambda_0$. If  $\alpha_i > 0$ then $(\lambda-\lambda_0)^{\alpha_i}$ is called an \emph{elementary divisor} at $\lambda_0$. Thus, to each $\lambda_0$, a set of increasingly ordered partial multiplicities $(\alpha_1, \ldots, \alpha_r)$ can be uniquely associated such that $0 \leq \alpha_1 \leq  \cdots \leq  \alpha_r$ jointly with a collection of elementary divisors $(\lambda-\lambda_0)^{\alpha_i}$ for $\alpha_i > 0$, including repetitions. The sum $\sum_{i=1}^r\alpha_i$ is the \emph{algebraic multiplicity} of the eigenvalue $\lambda_0$, while the number of nonzero terms in this sum is its \emph{geometric multiplicity}. An eigenvalue $\lambda_0$ is said to be \emph{simple}, if its algebraic multiplicity is one.

The sum of all partial multiplicities gives the total number of all finite eigenvalues of $P(\lambda)$ and is denoted as $\delta_{fin}(P)$. This value can be alternatively defined using the degrees of the invariant polynomials as follows
\[ \delta_{fin}(P) =  \sum_{i=1}^r \deg d_i(\lambda) . \]

For the definition of infinite eigenvalues of $P(\lambda)$ we use the mathematical framework introduced in \cite{Gohb09}, which we call the \emph{GLR framework} (using the initials of authors' names). For $j \geq d$, the $j$-\emph{reversal} of $P(\lambda)$ is the matrix polynomial   $\rev_j P(\lambda) := \lambda^jP(1/\lambda)$. If $j = d$, the $d$-reversal is called simply the \emph{reversa}l of $P(\lambda)$ and denoted $\rev P(\lambda)$. The GLR framework defines, for a grade $k$ polynomial matrix of degree $d$, $\lambda_0 = \infty $ an infinite eigenvalue of $P(\lambda)$ if and only if $0$ is an eigenvalue of $\rev_k P(\lambda)$. Using the Smith form of $\rev_k P(\lambda)$, we can define the increasingly ordered  partial multiplicities of the infinite eigenvalue as $(\alpha_1^\infty, \ldots, \alpha_r^\infty)$ with $0 \leq \alpha_i^\infty \leq \cdots \leq \alpha_r^\infty$ being the ordered partial multiplicities of the zero eigenvalue of $\rev_k P(\lambda)$. For each $\alpha_i^\infty > 0$ there exists an \emph{infinite elementary divisor} of degree $\alpha_i^\infty$ (or and infinite eigenvalue of multiplicity $\alpha_i^\infty$).  The number of infinite eigenvalues of $P(\lambda)$ is given by
\[ \delta_{\infty}(P) = \sum_{i=1}^r \alpha_i^\infty .\]
From the construction of the reversal follows that $P(\lambda)$ of grade $k = d$ has an eigenvalues at $\infty$ if and only if the rank of the leading coefficient matrix $P_d$ is strictly less than $r$. For a regular polynomial matrix this simply means that $P_d$ is singular. If $k > d$, then $P_k = 0$ and $P(\lambda)$ necessarily has infinite eigenvalues.
The following straightforward result (see \cite{DeTe14}) relates, in a simple way, the partial multiplicities of $P(\lambda)$ regarded as a grade $d$ polynomial matrix to the partial multiplicities of $P(\lambda)$ regarded as a grade $k$ polynomial matrix, with $k \geq d$.

\begin{lemma}\label{infshifts} Suppose $P(\lambda)$ is a polynomial matrix with rank $r$, degree $d$, grade $k = d$, and with the partial multiplicities $(\alpha_1^\infty, \ldots, \alpha_r^\infty)$ at $\infty$.  Then $P(\lambda)$ regarded as a polynomial matrix with grade $k \geq d$ has the partial multiplicities $\big(\alpha_1^\infty + (k-d), \ldots, \alpha_r^\infty+(k-d) \big)$  at $\infty$.
\end{lemma}

\section{Rational matrices} Let $R(\lambda) \in \mathds{F}(\lambda)^{m\times n} $ be an $m\times n$ \emph{rational matrix}  defined as
\be\label{ratmat} R(\lambda) = \ba{ccc} r_{11}(\lambda) & \cdots & r_{1n}(\lambda) \\
\vdots & \ddots & \vdots \\
r_{m1}(\lambda) & \cdots & r_{mn}(\lambda) \ea \, , \ee
where each $r_{ij}(\lambda)$ is a rational function (i.e., a ratio of two polynomials) of the form
\be\label{ratij} r_{ij}(\lambda) = \frac{\alpha_{ij}(\lambda)}{\beta_{ij}(\lambda)} = \frac{a_k\lambda^k+a_{k-1}\lambda^{k-1}+\cdots + a_1\lambda+a_0}{b_l\lambda^l+b_{l-1}\lambda^{l-1}+\cdots + b_1\lambda+b_0}  \ee
with coefficients in $\mathds{F}$. A polynomial matrix can be seen as a particular rational matrix with all denominator polynomials  equal to one (i.e., $\beta_{ij}(\lambda) = 1 \; \forall i, j$).

The following definitions are straightforward extensions of familiar notions for constant matrices.  $R(\lambda)$ is called \emph{invertible} if $m = n$ and $\det R(\lambda) \not\equiv 0$. Otherwise, $R(\lambda)$ is called \emph{singular}. The normal rank of $R(\lambda)$, denoted $\rank P(\lambda)$,  is the number of linearly independent rows or columns of $R(\lambda)$. A rational matrix $R(\lambda)$ is called \emph{proper} if $\lim_{\lambda \rightarrow \infty} R(\lambda) = D$, with $D$ having a finite norm. Otherwise, $R(\lambda)$ is called \emph{improper}. If $D = 0$, then $R(\lambda)$ is called \emph{strictly proper}. An invertible $R(\lambda)$ is \emph{biproper} if both $R(\lambda)$ and $R^{-1}(\lambda)$ are proper.

 For a rational matrix $R(\lambda)$ the so-called the McMillan framework, is widely used in control system theory to characterize the \emph{pole-zero structure} of $R(\lambda)$ \cite{Kail80}, \cite{Rose70}. In a broad sense, a complex value $\lambda_0$ is a pole of $P(\lambda)$ if at least one entry of $R(\lambda_0)$ is infinite, while $\lambda_0$ is a zero if $R(\lambda_0)$ has rank less than $r$ (its normal rank). This interpretation of poles and zeros leads to conceptual difficulties if $\lambda_0$ is both a pole and zero or if $\lambda_0 = \infty$ and therefore we give precise definitions based on the so-called local Smith-McMillan form (see, for example, \cite{Kail80}).

%If $\lambda_0$ is a finite value, then the pole-zero structure of $P(\lambda)$ at $\lambda = \lambda_0$, can be studied using the so-called local Smith-McMillan form.
\begin{theorem}[Local Smith-McMillan form at $\lambda_0$]\label{SMMF}
Let $R(\lambda)$ be an $m\times n$ rational matrix of rank $r$ with coefficients in  $\mathds{F}$ and $\lambda_0$ any finite value in $\overline{\mathds{F}}$. Then, there exist rational matrices $U_0(\lambda) \in \mathds{F}(\lambda)^{m\times m}$ and $V_0(\lambda) \in \mathds{F}(\lambda)^{n\times n}$, both regular at $\lambda_0$, such that
 \be\label{SmithMMForm} D_0(\lambda) := U_0(\lambda)R(\lambda)V_0(\lambda) =
 \ba{ccc|c} (\lambda-\lambda_0)^{\sigma_1} &&& \\ & \ddots && 0_{r,n-r} \\ &&(\lambda-\lambda_0)^{\sigma_r} & \\ \hline & 0_{m-r,r} && 0_{m-r,n-r} \ea , \ee
 where $\sigma_1 \leq \cdots \leq \sigma_r$. Moreover, $D_0(\lambda)$ is unique and is called the local Smith-McMillan form of $R(\lambda)$ at $\lambda_0$.
\end{theorem}
The values $\sigma_i$, $i = 1, \ldots, r$ are called the \emph{finite structural indices} at $\lambda_0$ and have the following interpretation.
A  value $\sigma_i<0$ defines a \emph{finite pole} of $R(\lambda)$ at $\lambda_0$ of multiplicity $-\sigma_i$, while a value $\sigma_i > 0$ defines a \emph{finite zero} of multiplicity $\sigma_i$ of $R(\lambda)$ at $\lambda_0$. $\lambda_0$ is neither pole nor zero if all structural indices are zero.
We denote $\delta^z_{fin}(R)$ the number of all finite zeros with their multiplicities, which is the  sum of all positive structural indices for $\lambda_0 \in \overline{\mathds{F}}$ and  denote $\delta^p_{fin}(R)$ the number of all finite poles, which is the absolute value of the sum of all negative structural indices for $\lambda_0 \in \overline{\mathds{F}}$.
%\[ \delta^z_{fin}(P) = \sum_{i=1}^r \min(\sigma_i,0) \]
%and
%\[ \delta^p_{fin}(P) = \sum_{i=1}^r \max(-\sigma_i,0) ,\]
%respectively.

 For a polynomial matrix $P(\lambda)$, all structural indices are non-negative, and therefore $P(\lambda)$ has no finite poles. It follows that $\delta^p_{fin}(P) = 0$. The following straightforward result states that the finite structural indices of a polynomial matrix $P(\lambda)$ are basically the same as the partial multiplicities of its finite eigenvalues.
\begin{lemma}\label{L:finzer} Let $P(\lambda)$ by a polynomial matrix of rank $r$ and let $\lambda_0$ be a finite eigenvalue  of $P(\lambda)$ with $(\alpha_1, \ldots, \alpha_r)$, the associated set of increasingly ordered partial multiplicities. Also,  let $(\sigma_1, \ldots, \sigma_r)$ be the set of increasingly ordered structural indices of $P(\lambda)$ at $\lambda_0$. Then, $\alpha_i = \sigma_i$ for $i = 1, \ldots, r$.
\end{lemma}

A similar result holds for the infinite poles and zeros.

\begin{theorem}[Local Smith-McMillan form at $\infty$]\label{SMMFinf}
Let $R(\lambda)$ be an $m\times n$ rational matrix of rank $r$ with coefficients in  $\mathds{F}$. Then, there exist rational matrices $U_\infty(\lambda) \in \mathds{F}(\lambda)^{m\times m}$ and $V_\infty(\lambda) \in \mathds{F}(\lambda)^{n\times n}$, both regular at $\infty$, such that
 \be\label{SmithMMForminf} D_\infty(\lambda) := U_\infty(\lambda)R(\lambda)V_\infty(\lambda) =
 \ba{ccc|c} (1/\lambda)^{\sigma_1^\infty} &&& \\ & \ddots && 0_{r,n-r} \\ &&(1/\lambda)^{\sigma_r^\infty} & \\ \hline & 0_{m-r,r} && 0_{m-r,n-r} \ea , \ee
 where $\sigma_1^\infty \leq \cdots \leq \sigma_r^\infty$. Moreover, $D_\infty(\lambda)$ is unique and is called the local Smith-McMillan form of $R(\lambda)$ at $\infty$.
\end{theorem}

The values $\sigma_i^\infty$, $i = 1, \ldots, r$ are called the \emph{infinite structural indices} and have a similar interpretation as before. A value $\sigma_i^\infty<0$ defines an infinite pole of $R(\lambda)$ of multiplicity $-\sigma_i^\infty$, while a value $\sigma_i^\infty > 0$ defines an infinite zero  of $R(\lambda)$ of multiplicity $\sigma_i^\infty$. $R(\lambda)$ has neither infinite poles nor infinite zeros if all infinite structural indices are zero.
We denote $\delta^z_{\infty}(R)$ the number of all infinite zeros with their multiplicities, which is the  sum of all positive infinite structural indices, and  denote $\delta^p_{\infty}(R)$ the number of all infinite poles, which is the absolute value of the sum of all negative infinite structural indices.

For a polynomial matrix $P(\lambda)$ all its poles are infinite, while its zeros may be both finite and infinite. The McMillan framework interprets  infinite zeros as ``infinite frequencies'' (e.g., as may occur in passive electrical networks), and therefore attaches a physically meaningful interpretation to infinite zeros. The following result shows that the relation between the infinite eigenvalues structure in the GLR framework and infinite zero structure in the McMillan framework can be expressed in term of a simple shift of multiplicities (see \cite{Ampa15}).

\begin{lemma}\label{L:infzer} Let $P(\lambda)$ by a polynomial matrix of rank $r$,  grade $k$ and let $(\alpha_1^\infty, \ldots, \alpha_r^\infty)$ be the set of increasingly ordered partial multiplicities associated to the infinite eigenvalues of $P(\lambda)$. Also,  let $(\sigma_1^\infty, \ldots, \sigma_r^\infty)$ be the set of increasingly ordered structural indices of $P(\lambda)$ at $\infty$. Then, $\sigma_i^\infty = \alpha_i^\infty-k$ for $i = 1, \ldots, r$.
\end{lemma}

If we know the partial multiplicities of the infinite eigenvalues of a polynomial matrix $P(\lambda)$, then we can simply determine the multiplicities of the infinite zeros from the positive structural indices $\sigma_j^{z,\infty} := \alpha^\infty_{j} - k$, $j = r-u+1, \ldots, r$, where $u$ is the number of partial multiplicities $\alpha^\infty_i$ which satisfy $\alpha^\infty_i > k$. In a similar way, we can determine the multiplicities of the infinite poles from the negative structural indices $\sigma_j^{p,\infty} := \alpha^\infty_{j}-k$, $j = 1, \ldots, l$, where $l$ is the number of partial multiplicities $\alpha^\infty_i$ which satisfy $\alpha^\infty_i < k$. Conversely, if we know the $\sigma_j^{p,\infty}$, $j = 1, \ldots, l$ and $\sigma_j^{z,\infty}$, $j = r-u+1, \ldots, r$, then for a grade $k$ polynomial matrix $P(\lambda)$, the partial multiplicities of infinite eigenvalues can be reconstructed as
\be\label{eigpolzer} (\sigma_1^{p,\infty}+k,\; \cdots, \; \sigma_l^{p,\infty}+k,\;  k,\;  \cdots , k,\;  \sigma_{r-u+1}^{z,\infty}+k, \; \cdots, \; \sigma_{r}^{z,\infty}+k ) ,  \ee
where there are $r-u-l$ partial multiplicities equal to $k$.
It must be noted that a consequence of \textsc{Lemma} \ref{infshifts} is, that, while the partial multiplicities of infinite eigenvalues depends on the chosen grade $k$ of the polynomial matrix $P(\lambda)$, the multiplicities of zeros and poles are independent of the choice of $k$. In particular, for a degree $d$ polynomial matrix, we always have $\sigma_1^{p,\infty} = -d$.

%If there are $k_\infty$ infinite eigenvalues of partial multiplicities $\alpha^\infty_i$ which satisfy $\alpha^\infty_i > d$ for $i = r-k_\infty+1, \ldots, r$, then there are $k_\infty$ infinite zeros of multiplicities $\sigma_j^{z,\infty} := \alpha^\infty_{r-k_\infty+j} - d$, $j = 1, \ldots, k_\infty$. Conversely, if there are $k_\infty$ infinite zeros of multiplicities $\sigma_j^{z,\infty}$, then the partial multiplicities of infinite eigenvalues contains $r-k_\infty$ values which satisfy $\alpha_j^\infty\leq d$, for $j = 1, \ldots, r-k_\infty$ and $k_\infty$ values which satisfy $ \alpha^\infty_{r-k_\infty+j} = \sigma_j^{z,\infty} + d$, for $j = 1, \ldots, k_\infty$.

%The number of infinite zeros of $P(\lambda)$  is given by
%\[  \delta_{\infty}^z(P) = \sum_{i=1}^r \max(\sigma_i^\infty,0) = \sum_{j=r-u+1}^{r} \sigma_j^{z,\infty} \]
%%and is related to the number of infinite eigenvalues as
%%\be\label{infeigzero}  \delta_{\infty}^z(P) = \delta_{\infty}(P)-k_\infty d - \delta_\infty^d(P), \ee
%%where $\delta_\infty^d(P)$ is the number of infinite eigenvalues of multiplicities at most $d$.
%%
%and the number of infinite poles of $P(\lambda)$ is given by
%\[  \delta_{\infty}^p(P) = \sum_{i=1}^r \max(-\sigma_i^\infty,0) = \sum_{j=1}^{l} \sigma_j^{z,\infty}  .\]
The number of finite and infinite poles of a rational matrix $R(\lambda)$, $\delta^p(R) := \delta_{fin}^p(R) + \delta_{\infty}^p(R)$, is called the \emph{McMillan degree} of $R(\lambda)$ \cite{Rose70} (also called the \emph{polar degree}). Analogously, the number of finite and infinite zeros is $\delta^z(R) := \delta_{fin}^z(R) + \delta_{\infty}^z(R)$ (also called the \emph{zero degree}).

%\begin{remark}
\emph{Remark. }Following the results of Verghese \cite{Verg80}, the pole structure of $R(\lambda)$ is equivalent to the zero structure of the regular rational matrix
\be\label{polzer} \widetilde P(\lambda) := \ba{cc} R(\lambda) & I_p \\ I_m & 0 \ea . \ee
Thus, we can convert the pole structure determination problem into a zero structure determination problem, which in turn can be solved as an eigenvalue computation problem. \finr
%\end{remark}

To characterize the singular structure of a rational matrix $R(\lambda)$, the relevant objects are the \emph{right nullspace} and \emph{left nullspace} of $R(\lambda)$.  For an $m\times n$ rational matrix $R(\lambda)$ of normal rank $r < \min(m,n)$,  consider the sets of left and right annihilators
\[ \mathcal{N}_l(R) := \{ v(\lambda) \in \mathds{F}(\lambda)^{1\times m}\; |\; v(\lambda)R(\lambda) = 0 \} , \]
\[ \mathcal{N}_r(R) := \{ v(\lambda) \in \mathds{F}(\lambda)^{n\times 1}\; |\; R(\lambda)v(\lambda) = 0 \} .\]
$\mathcal{N}_l(R)$ is a linear space of dimension $m-r$ called the
\emph{left nullspace} of $R(\lambda)$ and $\mathcal{N}_r(R)$
is a linear space of dimension $n-r$ called the
\emph{right nullspace} of $R(\lambda)$. It is always possible to choose polynomial bases $\{p_1(\lambda), \ldots, p_{m-r}(\lambda)\}$ and $\{q_1(\lambda), \ldots, q_{n-r}(\lambda)\}$ for $\mathcal{N}_l(R)$ and $\mathcal{N}_r(R)$, respectively. The degree of a polynomial basis is the sum of degrees of the basis polynomial vectors. A \emph{minimal polynomial basis} is one which has the least possible degree. For a minimal polynomial basis $\{p_1(\lambda), \ldots, p_{m-r}(\lambda)\}$ of the left nullspace $\mathcal{N}_l(R)$ the degrees $(\eta_1,\ldots,\eta_{m-r})$ of the polynomial vectors are called the \emph{left minimal indices} (also known as \emph{left Kronecker indices}), while for a minimal polynomial basis $\{q_1(\lambda), \ldots, q_{n-r}(\lambda)\}$ of the right nullspace $\mathcal{N}_r(R)$ the degrees $(\epsilon_1,\ldots, \epsilon_{n-r})$ of the polynomial vectors are called the \emph{right minimal indices} (also known as \emph{right Kronecker indices}). The left and right minimal indices are unique up to permutations and fully characterize the singular structure of a polynomial matrix. The above results have been established in \cite{Forn75} (see also \cite{Kail80} for a textbook presentation).

The degree of the minimal polynomial basis of $\mathcal{N}_l(R)$ is $\mu_l(R) := \sum_{i=1}^{m-r}\eta_{i}$ and, similarly, the degree of the minimal polynomial basis of $\mathcal{N}_r(R)$ is   $\mu_r(R):= \sum_{i=1}^{n-r}\epsilon_{i}$.  The sum of all the minimal indices of a given $R(\lambda)$ is
\[ \mu(R) := \mu_l(R) + \mu_r(R) . \]
If $R(\lambda)$ is invertible, then $\mu(R) := 0$. However, $\mu(R) := 0$ may generally occur (e.g., for a singular polynomial matrix (e.g., $R(\lambda) = A$, with $A$ singular).

There are several fundamental relationships between various structural elements of polynomial and rational matrices. The following result relates the finite and infinite eigenvalues of a regular pencil and is established in \cite[Lemma 6.1]{DeTe14}.

\begin{lemma} Let $P(\lambda)$ be  a regular $n\times n$ polynomial matrix of grade $k$, over an arbitrary field. Then
\[ \delta_{fin}(P) + \delta_\infty(P) = kn . \]
\end{lemma}

The above result is a corollary of the following more general relation involving the infinite eigenvalues, zeros and poles.
\begin{lemma}\label{L:V} Let $P(\lambda)$ be  an $m\times n$ polynomial matrix of grade $k$, rank $r$, over an arbitrary field. Then
\be\label{IST-GLRinf} \delta_{fin}(P) + \delta_{\infty}(P)  = kr + \delta^z(P)-\delta^p(P). \ee
\end{lemma}
For the  proof of this result we can apply the results of \textsc{Lemma} \ref{L:finzer} and  \textsc{Lemma}  \ref{L:infzer}, observing that
\[ \delta_{\infty}(P)  = kr + \delta^z_\infty(P)-\delta^p_\infty(P) \]
and taking into account that $\delta^p(P) = \delta^p_\infty(P)$.

The following result of \cite[Theorem 3]{Verg79} relates the number of  poles, number of zeros and the singular structure of a rational matrix.
\begin{lemma}\label{L:Verg} Let $R(\lambda)$ be  an $m\times n$ rational matrix over an arbitrary field. Then
\be\label{IST-VD} \delta^p(R) = \delta^z(R) + \mu(R) . \ee
\end{lemma}

The following result, called in \cite{DeTe14} the \emph{Index Sum Theorem},  relates the eigenvalue and singular structures of polynomial matrices.

\begin{lemma} Let $P(\lambda)$ be  an $m\times n$ polynomial matrix of grade $k$, rank $r$, over an arbitrary field. Then
\be\label{IST-GLR} \delta_{fin}(P) + \delta_{\infty}(P) + \mu(P) = kr . \ee
\end{lemma}
This result is Theorem 6.5 in \cite{DeTe14} and its proof is given in terms of companion form linearizations of the polynomial matrix $P(\lambda)$. An alternative, much simpler proof is possible by combining the results of \textsc{Lemma} \ref{L:V} and \textsc{Lemma} \ref{L:Verg}.

The handling of the particular case of a constant polynomial matrix $P(\lambda) := P_0$ depends on the choice of grade $k$. For $k = 0$, the polynomial matrix $P(\lambda)$ of rank $r = \rank \, P_0$ satisfies $P(\lambda) = \rev_0 P(\lambda)$ and therefore both $P(\lambda)$ and $\rev_0 P(\lambda)$ have the trivial Smith-form $\diag\{I_r,0\}$. It follows, that $P(\lambda)$ has no finite and infinite eigenvalues, and has $m-r$ right Kronecker indices equal to 0 and $n-r$ left Kronecker indices equal to 0 (both sets may be empty).
Regarded as a grade $k \geq 1$ polynomial matrix $P(\lambda) = P_0+\lambda 0 + \cdots + \lambda^k 0$, $P(\lambda)$ has no finite eigenvalues, but has $kr$ infinite eigenvalues with partial multiplicities $(k,k,\ldots,k)$, and the same left and right Kronecker indices as above.

\section{Matrix pencils}
A matrix pencil $M-\lambda N$ is a grade one polynomial matrix, whose structural properties can be numerically investigated using numerically reliable pencil manipulation algorithms. This allows to determine the structural properties of polynomial and rational matrices via linearization techniques.

In what follows, we assume $\mathds{F}$ is an algebraically closed field (e.g., $\mathds{F} = \mathds{C}$). The basic mathematical tool for matrix pencils is the \emph{Kronecker canonical form} (KCF) obtained using strict equivalence transformations, which exhibits both the eigenvalue structure as well as the singular structure of the pencil. Recall that two pencils $M-\lambda N$ and $\widetilde M-\lambda \widetilde N$ with $M, N, \widetilde M, \widetilde N \in \mathds{F}^{m\times n}$ are \emph{strictly equivalent} if there exist two invertible matrices $U \in \mathds{F}^{m\times m}$ and $V \in \mathds{F}^{n\times n}$  such that
\be\label{pencil-equiv}
U(M-\lambda N)V = \widetilde M -\lambda \widetilde N .\ee

For a general (singular) pencil, the strict equivalence leads to the KCF.
\begin{lemma}\label{L-KCF}
Let $M-\lambda N$ be an arbitrary pencil with $M, N \in \mathds{F}^{m\times n}$ and $\mathds{F}$  an algebraically closed field. Then, there exist invertible matrices $U \in \mathds{F}^{m\times m}$ and $V \in \mathds{F}^{n\times n}$ such that
\be\label{Kronecker} U(M-\lambda N)V = \ba{ccc} K_r(\lambda) \\ & K_{reg}(\lambda) \\&&K_l(\lambda)\ea ,\ee
where:
\begin{enumerate}
\item[1)] The full row rank pencil $K_r(\lambda)$ has the form
\[ K_r(\lambda) = \diag \big(L_{\epsilon_1}(\lambda), L_{\epsilon_2}(\lambda), \cdots, L_{\epsilon_{\nu_r}}(\lambda) \big) \, , \]
with  $L_{i}(\lambda)$  ($i \geq 0$) an $i\times (i+1)$  bidiagonal pencil of form
\be\label{Liblocks} L_i(\lambda) = \ba{cccc} -\lambda & 1 \\ & \ddots & \ddots \\ && -\lambda & 1 \ea  \, ; \ee
\item[2)] The regular pencil $K_{reg}(\lambda)$ is in the Weierstrass canonical form
\be\label{Weierstrass}  K_{reg}(\lambda) = \ba{cc}  J_f-\lambda I  \\ & I-\lambda  J_\infty \ea \, ,\ee
%with $\widetilde J_f$ in a (complex) Jordan canonical form as in (\ref{Jordan}) and
%with $\widetilde J_\infty$ in a nilpotent Jordan form as in (\ref{Jordan-null});
where $J_f$ is in the Jordan canonical form
\be\label{Jordan} J_f = \diag \left(J_{s_1}(\lambda_1), J_{s_2}(\lambda_2), \ldots, J_{s_k}(\lambda_k) \right) \, ,\ee
with  $J_{s_i}(\lambda_i)$ an elementary $s_i\times s_i$ Jordan block of the form
\[ J_{s_i}(\lambda_i) = \ba{cccc} \lambda_i & 1  \\ & \lambda_i & \ddots \\ & & \ddots & 1 \\ &&&\lambda_i \ea \]
and $J_\infty$ is nilpotent and has the (nilpotent) Jordan form
\be\label{Jordan-null} J_\infty = \diag \big(J_{s_1^\infty}(0), J_{s_2^\infty}(0), \ldots, J_{s_h^\infty}(0) \big) \, ;\ee
\item[3)] The full column rank $K_l(\lambda)$ has the form
\[ K_l(\lambda) = \diag \big(L^T_{\eta_1}(\lambda), L^T_{\eta_2}(\lambda), \cdots, L^T_{\eta_{\nu_l}}(\lambda) \big) \, .\]
\end{enumerate}
\end{lemma}

The Kronecker canonical form (\ref{Kronecker}) exhibits the right and left singular structures of the pencil $M-\lambda N$ via the full row rank block $K_r(\lambda)$ and full column rank block $K_l(\lambda)$, respectively, and the eigenvalue structure via the regular pencil $K_{reg}(\lambda)$.

The full row rank pencil $K_r(\lambda)$ is $n_r\times (n_r+\nu_r)$, where $n_r = \sum_{i=1}^{\nu_r} \epsilon_i$, the full column rank pencil $K_l(\lambda)$ is $(n_l+\nu_l)\times n_l$, where $n_l = \sum_{j=1}^{\nu_l} \eta_j$, while the regular pencil $K_{reg}(\lambda)$ is $n_{reg}\times n_{reg}$, with $n_{reg} = n_f+n_\infty$, where $n_f$ is the number of finite eigenvalues of $J_f-\lambda I$ and $n_\infty$ is the number of infinite eigenvalues of $I-\lambda J_\infty$. The $\epsilon_i\times (\epsilon_i+1)$ blocks $L_{\epsilon_i}(\lambda)$ with $\epsilon_i \geq 0$
are the right elementary Kronecker blocks, and $\epsilon_i$, for $i = 1, \ldots, \nu_r$, are called the \emph{right Kronecker indices}. The $(\eta_i+1)\times \eta_i$ blocks $L^T_{\eta_i}(\lambda)$ with $\eta_i \geq 0$
are the left elementary Kronecker blocks, and $\eta_i$, for $i = 1, \ldots, \nu_l$,  are called the \emph{left Kronecker indices}.

The Weierstrass canonical form (\ref{Weierstrass}) exhibits the finite and infinite eigenvalues of the pencil $M-\lambda N$. Each $s_i\times s_i$ Jordan block  $J_{s_i}(\lambda_i)$ corresponds to a finite elementary divisor $(\lambda-\lambda_i)^{s_i}$ and, by including all multiplicities, there are $n_f = \sum_{i=1}^{k}s_i$ \emph{finite eigenvalues}.  Each $s_i^\infty\times s_i^\infty$ nilpotent Jordan block  $J_{s_i^\infty}(0)$ corresponds to an infinite elementary divisor of order $s_i^\infty$ and there are $n_\infty = \sum_{i=1}^{h}s_i^\infty$ \emph{infinite eigenvalues}. Infinite eigenvalues with $s_i^\infty =1$ are called \emph{simple infinite eigenvalues}.
If $M-\lambda N$ is regular, then there are no left- and right-Kronecker structures and the Kronecker canonical form is simply the Weierstrass canonical form.

The normal rank $r$ of the pencil $M-\lambda N$ results as
\[  r := \rank (M-\lambda N) = n_r+n_f+n_\infty + n_l .\]
We can also express the rank $\ell$ of $N$ as
\[ \ell := \rank N = n_r+n_f+\rank J_\infty +n_l= n_r+n_f + \sum_{i=1}^{h}(s_i^\infty-1) +n_l =   r-h. \]
Assuming $s_1^\infty \leq s_2^\infty \leq \cdots \leq s_h^\infty$, then the $r$ partial multiplicities of the infinite eigenvalues are
\be\label{mulinfMN}
(\alpha_1^\infty, \alpha_2^\infty, \ldots, \alpha_r^\infty)=(0, \ldots, 0, s_1^\infty, \ldots , s_h^\infty ) , \ee
where the first $\ell = r-h$ partial multiplicities are equal to zero.

The pole-zero structure at $\infty$ of the pencil $M-\lambda N$ can be retrieved from the KCF using the result of \cite[Theorem 2]{Door83a}.

\begin{lemma}\label{L:VDD} Let $M-\lambda N$ be  an $m\times n$ linear matrix pencil of normal rank $r$ and let $\ell = \rank N$. Then, assuming $0 <  s_1^\infty \leq s_2^\infty \leq \cdots \leq  s_{h}^\infty $ are the ordered sizes of the nilpotent Jordan blocks of $J_\infty$, then the structural indices at $\infty$ of the pencil $M-\lambda N$ are determined by the KCF (\ref{Kronecker}) as follows:
\[ (\sigma_1^\infty, \sigma_2^\infty,\ldots,\sigma_r^\infty) = (-1, \ldots, -1, s_1^\infty-1, \ldots,  s_{h}^\infty -1 ) , \]
where there are $\ell$ structural indices equal to $-1$.
%\begin{enumerate}
%\item[1)] $\sigma_i^\infty = -1$ for $i = 1, \ldots, \ell$;
%\item[2)] $\sigma_i^\infty = \tilde s_i^\infty$
%\end{enumerate}
\end{lemma}
It follows that $M-\lambda N$ has $\ell$ poles at $\infty$, all of multiplicities equal to one, while the number of infinite zeros is $\sum_{i=1}^{h} (s_i^\infty-1) = n_\infty-h$.

The computation of the Kronecker-canonical form may involve the use of ill-conditioned transformations and, therefore, is potentially numerically unstable. Fortunately,  alternative so-called \emph{Kronecker-like forms} (KLFs), allow to obtain basically the same (or only a part of) structural information on the pencil $M-\lambda N$ by employing exclusively unitary transformations if $\mathds{F} = \mathds{C}$ (i.e., $U^*U = I$ and $V^*V = I$) or  orthogonal transformations  if $\mathds{F} = \mathds{R}$ (i.e., $U^TU = I$ and $V^TV = I$).

An arbitrary pencil $M-\lambda N$ can be reduced using orthogonal or unitary transformations $U$ and $V$ to the block-upper triangular form \cite{Door79a}
\be\label{KLF} U(M-\lambda N)V = \ba{cccc} M_r-\lambda N_r & \ast & \ast & \ast \\ 0 & M_{\infty}-\lambda N_{\infty} & \ast & \ast \\ 0 & 0 & M_f-\lambda N_f & \ast \\0&0&0 & M_l-\lambda N_l\ea ,\ee
where
\begin{enumerate}
\item[1)] $M_r-\lambda N_r$ has full row rank for all $\lambda  \in \mathds{F}$, has only a right nullspace, and contains information on the right Kronecker indices;
\item[2)] $M_{\infty}-\lambda N_{\infty}$ is regular and contains information on the infinite elementary divisors (i.e., the multiplicities of infinite eigenvalues);
\item[3)] $M_f-\lambda N_f$ is regular with $N_f$ invertible and contains the finite elementary divisors (i.e., the finite eigenvalues);
\item[4)] $M_l-\lambda N_l$ has full column rank for all $\lambda \in \mathds{F}$, has only a left nullspace, and contains information on the left Kronecker indices.
\end{enumerate}
The KLF (\ref{KLF}) can be obtained using numerically stable pencil reduction algorithms as proposed in \cite{Door79a}, \cite{Beel88}, \cite{Demm93}, \cite{Oara97}, which at the same time determine the left and right Kronecker indices and the infinite elementary divisors of $M-\lambda N$ from the fine block structure of subpencils $M_r-\lambda N_r$, $M_{\infty}-\lambda N_{\infty}$, and $M_l-\lambda N_l$. The finite eigenvalues can be computed using the QZ algorithm to compute the generalized eigenvalues of the pair $(M_f,N_f)$ \cite{Mole73}.

\emph{Remark.} The KLF (\ref{KLF}) separates the finite and infinite eigenvalues of $M-\lambda N$ as the eigenvalues of the regular subpencils $M_f-\lambda N_f$ and $M_\infty-\lambda N_\infty$, respectively, provides the information on the multiplicities of infinite eigenvalues (i.e, on the infinite elementary divisors of $M_\infty-\lambda N_\infty$), but does not provide further information on the multiplicities of the finite eigenvalues (i.e., on the finite elementary divisors of $M_f-\lambda N_f$).  For the determination of the partial multiplicities associated to a known finite eigenvalue $\lambda_0$ (e.g., computed using the QZ-algorithm), the following approach, suggested in \cite{Door79a}, can be employed. The pencil reduction algorithm is applied to the shifted pencil $N_f - \widetilde{\lambda} (M_f-\lambda_0 N_f)$ to determine its infinite elementary divisors. This corresponds to a transformation of the indeterminate as $\lambda = 1/(\widetilde\lambda -\lambda_0)$, which maps all finite eigenvalues at $\lambda_0$ of $M_f-\lambda N_f$ into infinite eigenvalues of $N_f -\widetilde{\lambda} (M_f-\lambda_0 N_f)$ for which the pencil reduction algorithm determines the partial multiplicities. \finr

The algorithms for the computation of Kronecker-like forms of linear pencils perform repeatedly column and row compressions of matrices using orthonal or unitary transformations. These operations involve rank determinations, for which rank revealing decompositions as the QR-decomposition with column pivoting or the more reliable (but also computationally more involved) singular value decomposition (SVD) can be used. The use of SVD-based rank determinations is the basis of the algorithms proposed in \cite{Door79a,Demm93}. Albeit numerically reliable, these algorithms have a computational complexity $\mathcal{O}(n^4)$, where $n$ is the minimum of row or column dimensions of the pencil.  More efficient algorithms of complexity $\mathcal{O}(n^3)$ have been proposed in \cite{Beel88,Oara97}, which rely on using QR decompositions with column pivoting for rank determinations. An enhanced version of algorithm of \cite{Oara97} can be devised by combining QR-decompositions (without column pivoting) and SVD-based rank determinations. Both compression techniques have been employed in the implementations of the basic tools to compute various KLFs in the \textsf{MatrixPencils} package along the lines of procedures described in \cite[see Procedure \textbf{PREDUCE}, Section 10.1.6]{Varg17}. Functions are also available for several applications of Kronecker-like forms as the computation of Kronecker indices, finite and infinite eigenvalues and zeros, normal rank. These functions served as building blocks for the implemented software for handling polynomial and rational matrices.

\section{Linearizations}
The standard way to address eigenvalue and structural analysis problems of matrix polynomials is via a \emph{linearization}, which replaces a given polynomial matrix $P(\lambda)$ by a matrix pencil $L(\lambda) = M - \lambda N$, which (ideally) preserves the eigenvalue  and singular structures of $P(\lambda)$. The structural analysis problems for $L(\lambda)$ are then solved using pencil reduction techniques in conjunction with the QZ-algorithm, as described in the previous section. Depending on the employed linearization, the structural properties of $P(\lambda)$ are retrieved from those of $L(\lambda)$.

Assume $P(\lambda)$ is a $p\times m$ polynomial matrix of grade $k$. A pencil $L(\lambda)$ is called a linearization of $P(\lambda)$ if there exist unimodular matrices $U(\lambda)$ and $V(\lambda)$ and $s \geq 0$ such that
\[ U(\lambda)L(\lambda)V(\lambda) = \diag \{ P(\lambda), I_s \} . \]
Thus, a linearization $L(\lambda)$ preserves the finite elementary divisors and thus the finite eigenvalues of $P(\lambda)$. It also preserves the dimensions of the right and left nullspaces of $P(\lambda)$. If in addition, $\rev_1 L(\lambda)$ is a linearization of $\rev_k P(\lambda)$, then $L(\lambda)$ is said to be a \emph{strong linearization} of $P(\lambda)$. For a strong linearization the infinite elementary divisors are also preserved. Therefore, the key property of a strong linearization is that $L(\lambda)$ and $P(\lambda)$ have the same finite and infinite elementary divisors. However, for a singular $P(\lambda)$ other structural features are also desirable to be preserved by $L(\lambda)$.

Among many existing strong linearizations, the Frobenius companian form linearizations are widely used in solving eigenvalue problems of polynomial matrices. A main appeal of these linearization is that they can be directly constructed from the coefficient matrices of $P(\lambda)$.  Besides the preservation of finite and infinite eigenvalue structures, these linearizations allow to easily retrieve information on the minimal indices. A potential drawback of these linearizations is that they usually do not reflect any structural feature which may be present in $P(\lambda)$ (e.g., symmetry).

A second category of linearizations is suitable for the investigation of the pole-zero and singular structures of a rational matrices using the McMillan framework. For a given rational matrix $R(\lambda)$, these linearizations are built as least order system matrix pencils of the form
\be\label{psysmat0} S(\lambda) = \ba{cc}-T(\lambda) & U(\lambda) \\ V(\lambda) & W(\lambda) \ea , \ee
where $T(\lambda)$, $U(\lambda)$, $V(\lambda)$ and $W(\lambda)$ are polynomial matrices of degree at most one,  with $T(\lambda)$ invertible, and satisfy the relation
\be\label{RosTFM0} R(\lambda) = V(\lambda)T^{-1}(\lambda)U(\lambda) + W(\lambda) . \ee
The notion of strong linearization can be extended to this framework, by requiring that $S(\lambda)$ preserves the complete pole-zero structure and singular structure of $R(\lambda)$. Two linearizations in this category are the pencil based linearization and the descriptor system based linearization. Important computational ingredients to determine these linearizations are minimal realization algorithms specific to each type of linearization.

\subsection{Companion forms based linearizations of polynomial matrices} These linearizations are widely used in the numerical linear algebra community, where the GLR framework is mostly employed. Assume the $p\times m$ polynomial matrix $P(\lambda)$ is given as a grade $k$ matrix polynomial of the form
\[ P(\lambda) = P_0 + P_1\lambda + \ldots + P_k \lambda^k .\]
The degree $d$ of $P(\lambda)$ is the maximum value of $i = 0, 1, \ldots, k$ for which $P_i \neq 0$ and $d \leq k$.

The \emph{first Frobenius companion form linearization} of $P(\lambda)$ is the linear pencil $C_1(\lambda) := M_1-\lambda N_1$, with
\be\label{CF1} M_1 = \ba{cccc} -P_{k-1} & -P_{k-2} & \cdots & -P_0 \\
 I_m & 0 & \cdots & 0 \\
  & \ddots & \ddots & \vdots \\
  0 && I_m & 0 \ea , \qquad N_1 = \ba{cccc} P_k \\ & I_m \\ && \ddots \\ &&& I_m \ea ,\ee
where $M_1$ and $N_1$ are $\big(p+(k-1)m\big)\times km$ matrices. If $P(\lambda)$ is regular then $C_1(\lambda)$ is regular as well. This linearization can be employed to recover the eigenvalue structure, zero structure, and singular structures of $P(\lambda)$ from the Kronecker structure of $C_1(\lambda)$ using the following results \cite{DeTe14}:
\begin{proposition}\label{PropCF1} Let $P(\lambda) =\sum_{i=0}^kP_i\lambda^i$ be a $p\times m$ matrix polynomial with grade $k \geq 2$,
and let $C_1(\lambda)$ be its first Frobenius companion form linearization. Then:
 \begin{enumerate}
\item[(a)] the finite and infinite elementary divisors of $P(\lambda)$ and $C_1(\lambda)$ are the same, thus $\delta_{fin}(P) = \delta_{fin}(C_1)$ and $\delta_{\infty}(P) = \delta_{\infty}(C_1)$;
%\item[(b)] if $r = \rank P(\lambda)$ and $r_1 = \rank C_1(\lambda)$, then the partial multiplicities $(\alpha_1^\infty, \ldots, \alpha_r^\infty)$ of the infinite eigenvalues of $P(\lambda)$ and the partial multiplicities $(\widetilde\alpha_1^\infty, \ldots, \widetilde\alpha_{r_1}^\infty)$ of the infinite eigenvalues of $C_1(\lambda)$  are related as
%    \[ \widetilde\alpha_i^\infty = 0, \; i = 1, \ldots, r_1-r, \quad \widetilde\alpha_{r_1-r+i}^\infty =  \alpha_i^\infty , \; i = 1, \ldots, r ;  \]
\item[(b)] if $r = \rank P(\lambda)$ and $r_1 = \rank C_1(\lambda)$, then the structural indices of $P(\lambda)$ at $\infty$ $(\sigma_1^{\infty}, \ldots, \alpha_{r}^{\infty})$ and the partial multiplicities $(\widetilde\alpha_1^\infty, \ldots, \widetilde\alpha_{r_1}^\infty)$ of the infinite eigenvalues of $C_1(\lambda)$  are related as
    \[ \widetilde\alpha_i^\infty = 0, \; i = 1, \ldots, r_1-r, \quad \widetilde\alpha_{r_1-r+i}^\infty =  \sigma_i^\infty+k , \; i = 1, \ldots, r ;  \]
\item[(c)] the right minimal indices $(\epsilon_1,\ldots, \epsilon_{\nu_r})\!$ of $P(\lambda)$ and the right minimal indices $(\widetilde\epsilon_1,\ldots, \widetilde\epsilon_{\nu_r})\!$ of $C_1(\lambda)$ are related as
    \[ \epsilon_i = \widetilde\epsilon_i - (k-1) , \; i = 1, \ldots, \nu_r ; \]
\item[(d)] the left minimal indices of $P(\lambda)$ and $C_1(\lambda)$ are the same, and hence
\[ \mu(P) = \mu(C_1)- (k-1)\nu_r; \]
\item[(e)] the normal ranks of $P(\lambda)$ and $C_1(\lambda)$ are related as
\[ \rank P(\lambda) = \rank C_1(\lambda) - m(k-1) .\]
 \end{enumerate}
 \end{proposition}

The \emph{second Frobenius companion form linearization} of $P(\lambda)$ is the linear pencil $C_2(\lambda) := M_2-\lambda N_2$, with
\be\label{CF2}  M_2 = \ba{cccc} -P_{k-1} & I_p &  & 0 \\
 -P_{k-2} & 0 & \ddots  \\
 \vdots & \vdots & \ddots & I_p \\
 -P_0 & 0 & \cdots & 0 \ea , \qquad N_2 = \ba{cccc} P_k \\ & I_p \\ && \ddots \\ &&& I_p \ea ,\ee
where $M_2$ and $N_2$ are $pk\times \big(m+(k-1)p\big)$ matrices.
This linearization can be employed to recover the eigenvalue structure, zero structure, and singular structures of $P(\lambda)$ from the Kronecker structure of $C_2(\lambda)$ using the following results \cite{DeTe14}:
\begin{proposition} Let $P(\lambda) =\sum_{i=0}^kP_i\lambda^i$ be a $p\times m$ matrix polynomial with grade $k \geq 2$,
and let $C_2(\lambda)$ be its second Frobenius companion form linearization. Then:
\begin{enumerate}
\item[(a)] the finite and infinite elementary divisors of $P(\lambda)$ and $C_2(\lambda)$ are the same, thus $\delta_{fin}(P) = \delta_{fin}(C_2)$ and $\delta_{\infty}(P) = \delta_{\infty}(C_2)$;
\item[(b)] if $r = \rank P(\lambda)$ and $r_2 = \rank C_2(\lambda)$, then the structural indices of $P(\lambda)$ at $\infty$ $(\sigma_1^{\infty}, \ldots, \alpha_{r}^{\infty})$ and the partial multiplicities $(\widetilde\alpha_1^\infty, \ldots, \widetilde\alpha_{r_2}^\infty)$ of the infinite eigenvalues of $C_2(\lambda)$  are related as
    \[ \widetilde\alpha_i^\infty = 0, \; i = 1, \ldots, r_2-r, \quad \widetilde\alpha_{r_2-r+i}^\infty =  \sigma_i^\infty+k , \; i = 1, \ldots, r ;  \]
\item[(c)] the right minimal indices of $P(\lambda)$ and $C_2(\lambda)$ are the same;
\item[(d)] the left minimal indices $(\eta_1,\ldots,\eta_{\nu_l})$ of $P(\lambda)$ and the left minimal indices $(\widetilde\eta_1,\ldots,\widetilde\eta_{\nu_l})$ of $C_2(\lambda)$ are related as
    \[ \eta_i = \widetilde\eta_i - (k-1) , \; i = 1, \ldots, \nu_l ; \]
and hence
\[ \mu(P) = \mu(C_2)-(k-1)\nu_l; \]
\item[(e)] the normal ranks of $P(\lambda)$ and $C_2(\lambda)$ are related as
\[ \rank P(\lambda) = \rank C_2(\lambda) - p(k-1) .\]
 \end{enumerate}
 \end{proposition}

\subsection{Linearizations of rational matrices}
To build a linearization of a rational matrix $R(\lambda)$ as a system matrix of the form (\ref{psysmat0}) we can use an additive decomposition of $R(\lambda)$ as
\be\label{addec0}
R(\lambda) = R_{sp}(\lambda) + R_{pol}(\lambda), \ee
where $R_{sp}(\lambda)$ is a strictly proper rational matrix and $R_{pol}(\lambda)$ is a polynomial matrix. As it will be shown below, we can build separate linearizations of the strictly proper part $R_{sp}(\lambda)$ and of the polynomial part $R_{pol}(\lambda)$, and combine them to obtain an overall linearization of $R(\lambda)$. We present the construction of two types of linearizations: the pencil based linearization and the descriptor system based linearization. Since the linearization of the strictly proper part is employed in both linearizations, we describe simple methods to build linearizations directly, by inspection, from the elements  of $R_{sp}(\lambda)$. For the linearization of the polynomial part, specific methods are described.
\subsubsection{Linearizations of strictly proper rational matrices}
For a strictly proper rational matrix $R_{sp}(\lambda)$, we show how to build a linearization, specified by the quadruple $(A_{sp} -\lambda I,B_{sp},C_{sp},0)$, which satisfies
\[ R_{sp}(\lambda) = C_{sp} \big(\lambda I -A_{sp}\big)^{-1} B_{sp} . \]
Let $r_{ij}(\lambda)$ denote the $(i,j)$-th entry of the strictly proper part $R_{sp}(\lambda)$ and let $n^{(ij)}$ be the degree of the denominator of $r_{ij}(\lambda)$. For each $r_{ij}(\lambda)$, a realization can be simply built (by inspection) in the form $(A^{(ij)} - \lambda I_{n^{(ij)}},b^{(ij)},c^{(ij)},0)$, which satisfies
\[ r_{ij}(\lambda) = c^{(ij)}\big(\lambda I_{n^{(ij)}} -A^{(ij)}\big)^{-1}b^{(ij)} . \]
Assume that $r_{ij}(\lambda)$ is expressed  as a ratio of two polynomials in the form
\be\label{glambda} r_{ij}(\lambda) =  \frac{\gamma_{k-1}\lambda^{k-1}+\cdots + \gamma_1\lambda+\gamma_0}{\lambda^k+\alpha_{k-1}\lambda^{k-1}+\cdots + \alpha_1\lambda+\alpha_0} \, ,\ee
where, for simplicity, we denoted $k := n^{(ij)}$. Then, $A^{(ij)}$, $b^{(ij)}$ and $c^{(ij)}$ can be expressed in a \emph{controllable} companion form
\[ \begin{array}{l} A^{(ij)} = \ba{ccccc} -\alpha_{k-1} & -\alpha_{k-2} & \cdots & -\alpha_1 & -\alpha_0 \\
1 & 0 & \cdots & 0 & 0\\
0 & 1 & \cdots & 0 & 0\\
\vdots & \vdots & \ddots & \vdots & \vdots \\
0 & 0 & \cdots & 1 & 0 \ea , \qquad b^{(ij)} = \ba{c} 1 \\0 \\ 0 \\ \vdots \\ 0 \ea , \\ \\
c^{(ij)} = \;\ba{ccccc} ~~\gamma_{k-1} & ~~\gamma_{k-2} & ~~\cdots & ~~\gamma_1 & ~~\gamma_0~ \ea .
\end{array}\]
This linearization has least order, provided the entry $r_{ij}(\lambda)$ in (\ref{glambda}) is in a cancelled form (i.e., with coprime numerator and denominator polynomials).
The matrices of the overall linearization are given by
\be\label{Rspnonmin} \begin{array}{l} A_{sp} = \ba{cccc} A^{(11)} \\ & A^{(12)} \\ &&  \ddots\\ &&& A^{(pm)} \ea, \qquad
B_{sp} = \ba{c} B^{(11)} \\ B^{(12)}\\  \vdots \\ B^{(pm)} \ea, \\ \\
C_{sp} = \ba{cccc} C^{(11)} & C^{(12)}&  \cdots & C^{(pm)} \ea ,
\end{array}\ee
where
\[ B^{(ij)} = [\; \underbrace{0 \; \cdots \; 0}_{j-1} \;\; b^{(ij)} \;\; \underbrace{0 \; \cdots \; 0}_{m-j} \;] , \qquad  C^{(ij)} = \ba{c} 0 \\ \vdots \\ 0 \\ c^{(ij)} \\ 0\\ \vdots \\ 0 \ea \hspace*{-6mm}
\begin{array}{l} \left.  \begin{array}{l} \\[-4mm] \\ \\ \\[-1mm] \end{array} \right\} \\  \\[2mm] \left.  \begin{array}{l} \\[-4mm] \\ \\ \\[-1mm] \end{array} \right\} \end{array} \hspace*{-3mm}\begin{array}{c}   \\[0mm] i-1 \\  \\ \\ \\ \\[-2mm] p-i \\ \\ \end{array} .\]
The overall linearization is generally not of least order (i.e., not controllable and not observable) and, therefore, to obtain a least order linearization of the strictly proper part an additional order reduction step is necessary. This step can be performed using numerically stable orthogonal transformation based algorithms, as those described in \cite{Door81}.

Using some preprocessing of the elements of $R_{sp}(\lambda)$, it is possible to construct linearizations which are either controllable or observable. In what follows, we only present a method to build a controllable linearization. For the computation of an observable linearization, a dual version of this method can be employed. Assume $R_{sp}(\lambda)$ is in a form with all entries in a column $j$ having the form $r_{ij}(\lambda) = h_{ij}(\lambda)/d_j(\lambda)$, with $d_j(\lambda)$  the least common multiple of the denominators of the elements of the $j$-th column. Then, for the $j$-th column of $R_{sp}(\lambda)$ we can build a linearization in the form $(\lambda I_{n^{(j)}} -A^{(j)},b^{(j)},C^{(j)},0)$, where $n^{(j)}$ is the degree of the common denominator of $d_{j}(\lambda)$. If we assume that $d_{j}(\lambda)$ is given by
\[ d_{j}(\lambda) = \lambda^k+\alpha_{k-1}\lambda^{k-1}+\cdots + \alpha_1\lambda+\alpha_0 \]
where $k := n^{(j)}$, and the numerators of the $j$-th columns can be expressed as a vector polynomial
\[ \ba{c} h_{1j}(\lambda) \\ \vdots \\ h_{pj}(\lambda) \ea =
\gamma_{k-1}\lambda^{k-1}+\cdots + \gamma_1\lambda+\gamma_0 ,\]
with each $\gamma_{i}$ a $p$-dimensional column vector, then $A^{(j)}$, $b^{(j)}$ and $C^{(ij)}$ can be expressed in a \emph{controllable} companion form
\[ \begin{array}{l} A^{(j)} = \ba{ccccc} -\alpha_{k-1} & -\alpha_{k-2} & \cdots & -\alpha_1 & -\alpha_0 \\
1 & 0 & \cdots & 0 & 0\\
0 & 1 & \cdots & 0 & 0\\
\vdots & \vdots & \ddots & \vdots & \vdots \\
0 & 0 & \cdots & 1 & 0 \ea , \qquad b^{(j)} = \ba{c} 1 \\0 \\ 0 \\ \vdots \\ 0 \ea , \\ \\
C^{(ij)} = \;\ba{ccccc} ~~\gamma_{k-1} & ~~\gamma_{k-2} & ~~\cdots & ~~\gamma_1 & ~~\gamma_0~ \ea .
\end{array}\]

The matrices of the overall linearization are given by
\be\label{spcontr} \begin{array}{l} A_{sp} = \ba{cccc} A^{(1)} \\ & A^{(2)} \\ &&  \ddots\\ &&& A^{(m)} \ea, \qquad
B_{sp} = \ba{c} B^{(1)} \\ B^{(2)}\\  \vdots \\ B^{(m)} \ea, \\ \\
C_{sp} = \ba{cccc} C^{(1)} & C^{(2)}&  \cdots & C^{(m)} \ea ,
\end{array}\ee
where
\[ B^{(j)} = [\; \underbrace{0 \; \cdots \; 0}_{j-1} \;\; b^{(j)} \;\; \underbrace{0 \; \cdots \; 0}_{m-j} \;] .\]
The overall linearization is controllable, but generally still not of least order (i.e., not observable). Therefore, to obtain a least order linearization of the strictly proper part an additional order reduction step is necessary, to remove the unobservable part. This step can be performed using numerically stable orthogonal transformation based algorithms, as those described in \cite{Door81}.
The dual version of this method builds an observable linearization and corresponds to apply the presented method to the transposed rational matrix $R_{sp}^T(\lambda)$. It may occasionally lead to lower dimensions and, therefore, may reduce the overall computational effort to determine a least order linearization.

\subsubsection{Pencils based linearization}
For a rational matrix $R(\lambda)$ (and therefore also for a polynomial matrix $P(\lambda)$) we can use a linearization of the form
\be\label{sysmat} S(\lambda) = \ba{cc} A- \lambda E & B- \lambda F \\ C - \lambda G & D - \lambda H \ea  , \ee
where $A-\lambda E$ is an $n\times n$ regular pencil and the quadruple of linear pencils $(A-\lambda E, B-\lambda F, C-\lambda G, D-\lambda H)$ is a \emph{pencils based realization} of $P(\lambda)$ which satisfies
\be\label{TFM1}      R(\lambda) = (C-\lambda G)(\lambda E - A)^{-1}(B-\lambda F) + D-\lambda H .\ee
In the control system literature, the matrix pencil $S(\lambda)$ is called the Rosenbrock's system matrix \cite{Rose70} of the pencils based realization of $R(\lambda)$.

Of particular interest are realizations which allow to retrieve the structural elements of $R(\lambda)$ from those of $S(\lambda)$. A realization $(A-\lambda E, B-\lambda F, C-\lambda G, D-\lambda H)$ is called \emph{strongly irreducible} \cite{Verg80} if it is \emph{strongly controllable}  and \emph{strongly observable}, for which the equivalent conditions are that the pencils
\be\label{strongmin} \ba{ccc} A-\lambda E & B-\lambda F & 0 \\ C-\lambda G & D-\lambda H & I_p \ea , \qquad \ba{cc} A-\lambda E & B-\lambda F \\ C-\lambda G & D-\lambda H \\ 0 & I_m \ea ,
\ee
have no finite and infinite zeros. These conditions are fulfilled if the pair $(A-\lambda E, B-\lambda F)$ is \emph{E-strongly controllable}  and the pair $(A-\lambda E, C-\lambda G)$  is \emph{E-strongly observable} for which the equivalent (stronger) conditions are that the pencils
\[ [\, A-\lambda E \;\; B-\lambda F ] , \qquad \ba{c} A-\lambda E \\ C-\lambda G \ea \]
have no finite and infinite eigenvalues \cite{Dopi20}.
In this case, the realization is called \emph{strongly minimal} and $n$  is the least achievable value such that (\ref{TFM1}) holds. This value is called in \cite{Rose70} the \emph{least order} and denoted with $\nu(P)$. The main importance of \emph{strongly minimal} realizations is that they can be computed in a relatively simple way from a non-minimal realization (e.g., using a procedure proposed in \cite{Dopi20}) using standard pencil manipulation algorithms.

To build a pencil based linearization of a rational matrix $R(\lambda)$ we additively  decompose $R(\lambda)$ as $R(\lambda) = R_{sp}(\lambda) + R_{pol}(\lambda)$, where
$R_{sp}(\lambda)$ is a strictly proper rational matrix and $R_{pol}(\lambda)$ is a polynomial matrix. As we have already shown, we can build for the strictly proper part $R_{sp}(\lambda)$ a linearization of the form  $(A_{sp} - \lambda I,B_{sp},C_{sp},0)$, while for the polynomial part $R_{p}(\lambda)$ we will show that we can build a pencil based linearization of the form $(A_{pol} - \lambda E_{pol},B_{pol}- \lambda F_{pol},C_{pol}- \lambda G_{pol},D_{pol}- \lambda H_{pol})$, such that the overall realization of $R(\lambda)$ is obtained in the form
\be\label{lplin1} R(\lambda) = \ba{c|c} A-\lambda E & B-\lambda F \\ \hline C-\lambda G & D-\lambda H \ea := \ba{cc|c} A_{sp}-\lambda I & 0 & B_{sp} \\ 0 & A_{pol}-\lambda E_{pol} & B_{pol}- \lambda F_{pol} \\ \hline
C_{sp} & C_{pol}- \lambda G_{pol} & D_{pol}- \lambda H_{pol} \ea \, .\ee

For the polynomial part $R_{pol}(\lambda)$, a pencil based linearization can be easily derived (by inspection) assuming a grade $k \geq 2$ polynomial matrix $R_{pol}(\lambda) = \sum_{i = 0}^k \lambda^i P_i$ as follows. A strongly-controllable realization of order $m(k-1)$ is given by
\be\label{scon} {\arraycolsep=1mm\ba{c|c} A_{pol}-\lambda E_{pol} & B_{pol}-\lambda F_{pol} \\ \hline C_{pol} - \lambda G_{pol} & D_{pol} - \lambda H_{pol} \ea  :=
\ba{ccccc|c} I_m & -\lambda I_m &&&& 0 \\ & I_m & -\lambda I_m &&& 0\\ && \ddots & \ddots && \vdots\\ &&&\ddots & -\lambda I_m & 0 \\ &&&& I_m & -\lambda I_m \\ \hline
 P_{k-1}+\lambda P_{k} & P_{k-2} && \cdots &  P_1 & P_0 \ea }.
\ee
A strongly-observable realization of order $p(k-1)$ is given by
\be\label{sobs} {\arraycolsep=1mm\ba{c|c} A_{pol}-\lambda E_{pol} & B_{pol}-\lambda F_{pol} \\ \hline C_{pol} - \lambda G_{pol} & D_{pol} - \lambda H_{pol} \ea  :=
\ba{ccccc|c} I_p & -\lambda I_p &&&& P_1  \\ & I_p & -\lambda I_p &&& P_2 \\ && \ddots & \ddots && \vdots\\ &&&\ddots & -\lambda I_p & P_{k-2} \\ &&&& I_p & P_{k-1}+\lambda P_{k} \\ \hline
 -\lambda I_p & 0 && \cdots &  0 & P_0 \ea }.
\ee
For $k = 1$, a realization of order $0$ is given by $D_{pol}-\lambda H_{pol} := P_0+\lambda P_1$, while for a constant polynomial matrix (i.e., $k = 0$), we take $D_{pol} := P_0$ and $H_{pol}$ an empty matrix.

The following result has been stated in \cite{Verg80} for strongly irreducible realizations (and thus also valid for strongly minimal realizations):
\begin{proposition}\label{Plp} Let $R(\lambda)$ be a $p\times m$ rational matrix
and let $(A-\lambda E, B-\lambda F, C-\lambda G, D-\lambda H)$ be a \emph{strongly irreducible} linearization satisfying (\ref{TFM1}). Then:
\begin{enumerate}
\item[(a)] the finite and infinite zero structures and the singular Kronecker structures of $R(\lambda)$ and $S(\lambda)$ are the same;
\item[(b)] the finite and infinite pole structures of $R(\lambda)$ and the finite and infinite zero structures of the \emph{pole pencil}
\be\label{TFM-poles} \widetilde S(\lambda) := \ba{ccc} A & B & 0\\ C & D & I_p \\ 0 & I_m & 0\ea - \lambda  \ba{ccc} E & F & 0 \\ G & H & 0 \\ 0 & 0 & 0 \ea \ee
are the same.
\end{enumerate}
\end{proposition}

\emph{Remark. }For computational purposes, the reduced pole pencil
\[ \widehat S(\lambda) = \ba{ccc} A & 0 & 0\\ 0 & 0 & I_p \\ 0 & I_m & 0\ea - \lambda  \ba{ccc} E & F & 0 \\ G & H & 0 \\ 0 & 0 & 0 \ea \]
can be employed instead of (\ref{TFM-poles}) to compute both the finite and infinite poles of $R(\lambda)$, as the zeros of  $\widetilde S(\lambda)$ using pencil manipulation techniques. The finite poles can alternatively be determined as the finite eigenvalues of $A-\lambda E$. If the underlying linearization of $R(\lambda)$ is of  the form (\ref{lplin1}), with $A_{sp}$ of least possible order, then the finite poles of $R(\lambda)$ are simply the eigenvalues of $A_{sp}$ and the infinite poles are the infinite zeros of the pencil
\be\label{Shat} \widehat S_\infty(\lambda) := \ba{ccc} A_{pol} & 0 & 0\\ 0 & 0 & I_p \\ 0 & I_m & 0\ea - \lambda  \ba{ccc} E_{pol} & F_{pol} & 0 \\ G_{pol} & H_{pol} & 0 \\ 0 & 0 & 0 \ea .\ee
If $R(\lambda)$ is a polynomial matrix, it has no finite poles. In this case, $E$ is nilpotent and $A-\lambda E$ is unimodular.  \finr

For a polynomial matrix $R(\lambda) = P(\lambda)$, with $P(\lambda)$ as given in (\ref{matpol}),  it is straightforward to relate the number of infinite eigenvalues of $P(\lambda)$ and of the associated system matrix $S(\lambda)$. Using (\ref{IST-GLR}), we have
 \[ \delta_{\infty}(P) - \delta_{\infty}(S) = dr-\rank(S). \]
 Note that for an $n$-th order realization, $\rank(S) = r+n$ and we have
 \be\label{stronginf} \delta_{\infty}(P) - \delta_{\infty}(S) = (d-1)r-n. \ee
It follows that, knowing $\delta_{\infty}(S)$, the number of infinite eigenvalues $\delta_{\infty}(P)$ can be recovered using (\ref{stronginf}). Alternatively, knowing the pole-zero structure (i.e., the multiplicities of infinite poles $\sigma_j^{p,\infty}$, $j = 1, \ldots, l$ and the multiplicities of infinite zeros $\sigma_j^{z,\infty}$, $j = r-u+1, \ldots, r$), then for a grade $k$ polynomial matrix $P(\lambda)$, the partial multiplicities of infinite eigenvalues can be reconstructed from (\ref{eigpolzer}).

Building linearizations of the form (\ref{sysmat}) for a rational matrix $R(\lambda)$ based on a strongly irreducible realization $(A-\lambda E, B-\lambda F, C-\lambda G, D-\lambda H)$ usually involves two stages. First, we need to determine a least order linearization of the strictly proper part of $R(\lambda)$ starting with an initial realization $(A_{sp} -\lambda I,B_{sp},C_{sp},0)$, which may be uncontrollable, or unobservable, or both uncontrollable and unobservable. The computation of the least order linearization can be performed using numerically stable orthogonal transformation based algorithms, as those described in \cite{Door81}. In the second stage, we  can build a strongly controllable realization of the polynomial part $R(\lambda)$ as in (\ref{scon}), which however may not be strongly observable, because the pencil $\left[ \begin{smallmatrix} A_{pol}-\lambda E_{pol}\\ C_{pol}-\lambda G_{pol} \end{smallmatrix} \right]$ may have infinite eigenvalues. These infinite eigenvalues can be removed using the procedure proposed in \cite{Dopi20}. A completely similar approach can be devised by starting with a strongly observable realization as in (\ref{sobs}) and then removing the infinite eigenvalues of the pencil $[\, A_{pol}-\lambda E_{pol} \; \; B_{pol}-\lambda F_{pol} \,]$ using the procedure of \cite{Dopi20}. The decision on which of these approaches to be used can be guided by the goal to minimize the computational effort in the second stage, by choosing the initial realization of lower order. Therefore, if $p > m$, the realization (\ref{scon}) of order $m(k-1)$ is to be preferred, while if $p < m$ the realization (\ref{sobs}) of order $p(k-1)$ may be preferable.

\subsubsection{Descriptor system based linearization}
For a rational matrix $R(\lambda)$ (and therefore also for a polynomial matrix $P(\lambda)$) we can alternatively use a linearization with a system matrix of the form
\be\label{sysmatD} S(\lambda) = \ba{cc} A -\lambda E & B \\ C & D \ea , \ee
where $(A-\lambda E, B, C, D)$, called a \emph{descriptor system realization} of $R(\lambda)$, satisfies
\be\label{TFM2}      R(\lambda) = C(\lambda E - A)^{-1}B + D .\ee
 The descriptor realization is called \emph{irreducible} if it is \emph{controllable} and  \emph{observable}, for which equivalent conditions are  that the pencils
 \[ [\, A-\lambda E \;\; B\,] , \qquad \ba{c} A-\lambda E \\ C \ea \]
 have no  finite and infinite zeros \cite{Verg81}. If additionally the pencil $A-\lambda E$ has no first order infinite elementary divisors (also called  \emph{non-dynamic modes}), then the descriptor realization is called \emph{minimal} and $n$, the order of $A$, is the least achievable dimension.

Any descriptor realization of $R(\lambda)$ is a particular pencil realization, which is strongly irreducible if the descriptor realization is irreducible. Therefore, the results of Proposition \ref{Plp} apply also to an irreducible descriptor system realization with the system matrix $S(\lambda)$ in (\ref{sysmatD}). An advantage of using descriptor system realizations for pole computations is that the finite and infinite poles of $R(\lambda)$ can be determined as the finite and infinite zeros of the reduced pole pencil
\be\label{redpol} \widetilde S(\lambda) := A-\lambda E . \ee

For a polynomial matrix $P(\lambda)$, using (\ref{stronginf}), the number of infinite eigenvalues can be recovered from those of the system matrix $S(\lambda)$ in (\ref{sysmatD}), while the partial multiplicities of infinite eigenvalues of $P(\lambda)$ can be retrieved from the infinite zero structures of $S(\lambda)$ and $\widetilde S_p(\lambda)$ in (\ref{redpol}).

To build a descriptor system based linearization of a rational matrix $R(\lambda)$ we additively  decompose $R(\lambda)$ as $R(\lambda) = R_{sp}(\lambda) + R_{pol}(\lambda)$, where
$R_{sp}(\lambda)$ is a strictly proper rational matrix and $R_{pol}(\lambda)$ is a polynomial matrix. As we have already shown, we can build for the strictly proper part $R_{sp}(\lambda)$ a linearization of the form  $(A_{sp} - \lambda I,B_{sp},C_{sp},0)$, while for the polynomial part $R_{p}(\lambda)$ we will show that we can build a linearization of the form $(I - \lambda E_{pol},B_{pol},C_{pol},D_{pol})$, such that the overall realization of $R(\lambda)$ is obtained in the form
\be\label{Rtot} R(\lambda) = \ba{c|c} A-\lambda E & B \\ \hline C & D \ea := \ba{cc|c} A_{sp}-\lambda I & 0 & B_{sp} \\ 0 & I-\lambda E_{pol} & B_{pol} \\ \hline
C_{sp} & C_{pol} & D_{pol} \ea \, .\ee

For the polynomial matrix $R_{pol}(\lambda)$ it is always possible to build a strongly irreducible realization of the form $(I - \lambda E_{pol},B_{pol},C_{pol},D_{pol})$,
with   $E_{pol}$ nilpotent (and therefore $I-\lambda E_{pol}$ unimodular). Such a realization can be determined following the suggestions from \cite{Verg79} by building a (standard) minimal realization $(E_{pol} -\lambda I,B_{pol},C_{pol},D_{pol})$ of the strictly proper rational matrix $\lambda^{-1}(R_{pol}(\lambda^{-1})-P_0)$ satisfying $\lambda^{-1}(R_{pol}(\lambda^{-1})-P_0) = -C_{pol}(\lambda I-E_{pol})^{-1}B_{pol}$. For this purpose, minimal realization procedures as suggested in \cite{Door81} can be used. Then, the realization of $P(\lambda)$ is simply $(I-\lambda E_{pol}, B_{pol}, C_{pol}, P_0)$.

Building linearizations of the form  $(I-\lambda E_{pol}, B_{pol}, C_{pol}, D_{pol})$  for the polynomial matrix $R_{pol}(\lambda) = \sum_{i = 0}^k \lambda^i P_i$   can be done in two steps. First, we build a controllable realization of $R_{pol}(\lambda)$ in the form
\be\label{sscon}  \ba{c|c} I - \lambda E_{pol} & B_{pol} \\ \hline C_{pol} & D_{pol} \ea := \ba{ccccc|c} I_m & -\lambda I_m &&&& 0 \\ & I_m & -\lambda I_m &&& 0\\ && \ddots & \ddots && \vdots\\ &&&I_m & -\lambda I_m & 0 \\ &&&& I_m & -I_m \\ \hline
 P_{k} & P_{k-1} & \cdots &  P_1 & 0 & P_0 \ea .
\ee
of order $m(k+1)$,
which however may not be observable at infinity, because the pair $\left[ \begin{smallmatrix} I-\lambda E_{pol}\\ C_{pol} \end{smallmatrix} \right]$ may have infinite (decoupling) zeros, or equivalently, the standard pair $(E_{pol},C_{pol})$ may have unobservable null eigenvalues. These unobservable eigenvalues can be removed by reducing the pair $(E_{pol},C_{pol})$ to the observability staircase form \cite{Door81} from which an observable realization can be obtained. A suitable algorithm for this purpose is described, for example, in \cite{Door81}.

A completely similar approach can be devised by starting at the first step with an observable realization of $R_{pol}(\lambda)$ in the form
\be\label{ssobs}  \ba{c|c} I - \lambda E_{pol} & B_{pol} \\ \hline C_{pol} & D_{pol} \ea := \ba{ccccc|c} I_p & -\lambda I_p &&&& 0 \\ & I_p & -\lambda I_p &&& P_1\\ && \ddots & \ddots && \vdots\\ &&&I_p & -\lambda I_p & P_{k-1} \\ &&&& I_p & P_{k} \\ \hline
 -I_p & 0 & \cdots &  0 & 0 & P_0 \ea .
\ee
of order $p(k+1)$, and then removing the uncontrollable infinite eigenvalues of the pencil $[\, I\!-\!\lambda E_{pol} \;  B_{pol} \,]$ using the procedure of \cite{Door81}. The decision on which of these approaches to be used can be guided by the goal to minimize the computational effort in the second step, by choosing the initial realization of lower order. Therefore, if $p > m$, the realization (\ref{sscon}) of order $m(k+1)$ is to be preferred, while if $p < m$ the realization (\ref{ssobs}) of order $p(k+1)$ may be preferable.

\subsection{Linearization of polynomial system matrices}
A polynomial system matrix has the form
\be\label{psysmat} S(\lambda) = \ba{cc}-T(\lambda) & U(\lambda) \\ V(\lambda) & W(\lambda) \ea , \ee
where $T(\lambda)$, $U(\lambda)$, $V(\lambda)$ and $W(\lambda)$ are polynomial matrices of sizes $n\times n$, $n\times m$, $p\times n$ and $p\times m$, respectively, and $T(\lambda)$ is regular. This polynomial system matrix is associated to the rational \emph{transfer function matrix}
\be\label{RosTFM} R(\lambda) = V(\lambda)T^{-1}(\lambda)U(\lambda) + W(\lambda)  \ee
and it was used in the works of Rosenbrock \cite{Rose70} and later of  Verghese \cite{Verg79,Verg80} to study the pole-zero and singular structures of $R(\lambda)$. Particular system matrices with first order polynomial matrices (i.e., pencils) have already been considered in (\ref{sysmat}) and (\ref{sysmatD}), and are well-suited for computational purposes. Therefore, the linearization of a general polynomial matrix to obtain a first order polynomial form is often necessary.

The linearization of the  polynomial system matrix $S(\lambda)$ in (\ref{psysmat}) to a first order polynomial matrix of the form  (\ref{sysmat}) having the same transfer function matrix $R(\lambda)$ can be performed using the pencil based linearizations of $S(\lambda)$ as in (\ref{scon}) or (\ref{sobs}). Assume that $S(\lambda)$ has the following pencil based linearization
\[ \widetilde C_1(\lambda) = {\def\arraystretch{1.4}\ba{cc|c} \widetilde A- \lambda \widetilde E & \widetilde B_1- \lambda \widetilde F_1 & \widetilde B_2- \lambda \widetilde F_2 \\ \widetilde C_1 - \lambda \widetilde G_1 & \widetilde D_{11} - \lambda \widetilde H_{11} & \widetilde D_{12} - \lambda \widetilde H_{12} \\ \hline
\widetilde C_2 - \lambda \widetilde G_2 & \widetilde D_{21} - \lambda \widetilde H_{21} & \widetilde D_{22} - \lambda \widetilde H_{22} \ea :=
\ba{c|c}  A- \lambda E & B- \lambda F \\ \hline
  C- \lambda G & D- \lambda H \ea } , \]
where $\widetilde A- \lambda \widetilde E$ is a regular $\tilde n\times \tilde n$ pencil with $\tilde n$ depending on the chosen linearization (\ref{scon}) or (\ref{sobs}),  $\widetilde D_{11} - \lambda \widetilde H_{11}$ is an $n\times n$ pencil, and $\widetilde D_{22} - \lambda \widetilde H_{22}$ is a $p\times m$ pencil (the dimensions of the rest of subpencils implicitly result). The resulting $(n+\tilde n)\times (n+\tilde n)$ pencil $A- \lambda E$ is regular and it can be shown that
\[ R(\lambda) = (C- \lambda G)(\lambda E -A)^{-1}(B- \lambda F)+ D- \lambda H .\]
This realization is usually not strongly minimal and the reduction to a least order linearization can be achieved using the techniques described in \cite{Dopi20}.
A similar approach can be devised to arrive to a descriptor system based linearization of the form (\ref{sysmatD}).

Several particular cases of interest of the transfer function matrix (\ref{RosTFM}) can be addressed in a straightforward manner by suitably defining the quadruple  $\{ T(\lambda), U(\lambda), V(\lambda), W(\lambda)\}$ of polynomial matrices which enter in (\ref{psysmat}). The \emph{left polynomial matrix fractional description }
\[ R(\lambda) = D^{-1}(\lambda)N(\lambda) \]
corresponds to the quadruple $\{ D(\lambda), N(\lambda), I, 0)\}$, while the \emph{right polynomial matrix fractional description }
\[ R(\lambda) = N(\lambda)D^{-1}(\lambda) \]
corresponds to the quadruple $\{ D(\lambda), I, N(\lambda),  0)\}$. Finally, for the inverse of a square polynomial matrix $P(\lambda)$
\[ R(\lambda) = P^{-1}(\lambda) , \]
the quadruple $\{ P(\lambda), I, I,  0)\}$ may serve as basis of constructing suitable linearizations.
\section{Implemented software} \label{app:software}

In what follows, we succinctly describe the newly implemented software tools for the releases \textsf{v1.0} and \textsf{v1.1} of the \textsf{Julia} package \textsf{MatrixPencils} \cite{Varg20}. These functions cover the computation of structural elements of polynomial and rational matrices and related computations as described in this paper. The required basic computational tools, as for example, tools for the computation of the Kronecker structure of matrix pencils or the computation of least order linearizations have been already implemented for the (previous) \textsf{Release v0.5} and will be described elsewhere. The implemented functions focus only on the
computation of structural elements such as eigenvalues, zeros, Kronecker indices, but do not address the computation of vectors associated with them, as eigenvectors, zero directions, or bases vectors of certain nullspaces.

A polynomial matrix $P(\lambda)$ can be entered in the \textsf{Julia} language in  two formats. The first possibility is to enter it as a matrix with \textsf{Polynomial} type elements as defined in the  \textsf{Polynomials} package (\url{https://github.com/JuliaMath/Polynomials.jl}). This input format is mainly intended as a convenient way to enter polynomial matrices in a quasi-symbolic form using matrices or vectors with polynomial entries (or even scalar polynomials) as input data.
The second input format relies on the monomial basis representation $P(\lambda) = P_0 + P_1\lambda + \ldots + P_k \lambda^k$  by storing the coefficient matrices $P_0$, $P_1$, $\ldots$ , $P_k$ of the successive powers $\lambda^0$, $\lambda^1$, $\ldots$, $\lambda^k$ in an 3-dimensional array \texttt{P}, where \texttt{P[:,:,i]} contains $P_{i+1}$. This format is internally used in all computational routines and is also suited for alternative representations of $P(\lambda)$ (e.g., in other polynomial bases).

A rational matrix $R(\lambda)$ can be entered using a pair of polynomial matrices $N(\lambda)$ and $D(\lambda)$ containing, respectively, the numerator and denominator polynomials of the elements of $R(\lambda)$. All elements of the denominator polynomial matrix $D(\lambda)$ must be nonzero. By convention, a polynomial matrix, viewed as a rational matrix, has all denominators equal to 1, and therefore only the numerator polynomial matrix $N(\lambda)$  has to be entered. Rational matrices can also be implicitly defined via polynomial system matrices as in (\ref{psysmat}), left or right polynomial matrix fraction descriptions $D^{-1}(\lambda)N(\lambda)$ and $N(\lambda)D^{-1}(\lambda)$, respectively, or as the inverse of a polynomial matrix $P(\lambda)$.

The following mnemonics have been used in the naming of functions:
\begin{center}
\begin{longtable}{ll} Mnemonic & Denotation \\ \hline
\texttt{\bfseries lp}     & linear pencil\\
\texttt{\bfseries ls}     & linear system in descriptor form\\
\texttt{\bfseries lps}     & linear pencil system\\
\texttt{\bfseries pm}     & polynomial matrix\\
\texttt{\bfseries rm}     & rational matrix\\
\texttt{\bfseries spm}     & structured polynomial matrix; also polynomial system matrix\\
\texttt{\bfseries lpmfd}     & left polynomial matrix fractional description\\
\texttt{\bfseries rpmfd}     & right polynomial matrix fractional description\\
\texttt{\bfseries poly}     & polynomial matrix, polynomial vector or scalar polynomial\footnote{Based on the \textsf{Polynomial} type
provided by the  \textsf{Polynomials} package \url{https://github.com/JuliaMath/Polynomials.jl}}\\
\texttt{\bfseries 2}     & place holder for ``conversion to''\\
\end{longtable}
\end{center}

The following table lists the main functions available for polynomial matrices in \textsf{Release v1.0} (and later) of the \textsf{MatrixPencils} package:
%\begin{table}[h]
\begin{longtable}{lp{10.5cm}} Function & Description \\ \hline
  \texttt{\bfseries poly2pm}     & Conversion of a polynomial matrix used in \textsf{Polynomials} package to
  a polynomial matrix represented as a 3-dimensional matrix\\
  \texttt{\bfseries pm2poly}     & Conversion of a polynomial matrix represented as a 3-dimensional matrix to a polynomial matrix used in \textsf{Polynomials} package\\
  \texttt{\bfseries pmdeg}     & Determination of the degree of a polynomial matrix \\
  \texttt{\bfseries pmeval}     & Evaluation of a polynomial matrix for a given value of its argument\\
  \texttt{\bfseries pmreverse}       & Building the reversal of a polynomial matrix  \\
  \texttt{\bfseries pmdivrem}\footnote{Included from \textsf{Release v1.1}}     & Evaluation of the quotients and remainders of the divisions of the numerators by the corresponding denominators of a rational matrix \\
  \texttt{\bfseries pm2lpCF1}    & Building a linearization in the first Frobenius companion form \\
  \texttt{\bfseries pm2lpCF2}  & Building a linearization in the second Frobenius companion form  \\
  \texttt{\bfseries pm2ls}    & Building a structured linearization $\left[ \begin{smallmatrix}A-\lambda E & B\\ C & D\end{smallmatrix}\right]$ of a polynomial matrix \\
  \texttt{\bfseries ls2pm}     & Computation of the polynomial matrix from its structured linearization\\
  \texttt{\bfseries pm2lps}      & Building a pencil based structured linearization
  $\left[ \begin{smallmatrix}A-\lambda E & B-\lambda F\\ C-\lambda G  & D-\lambda H\end{smallmatrix}\right]$
   of a polynomial matrix\\
  \texttt{\bfseries lps2pm}      & Computation of the polynomial matrix from its pencil based structured linearization  \\
  \texttt{\bfseries spm2ls}    & Building a structured linearization $\left[ \begin{smallmatrix}A-\lambda E & B\\ C & D\end{smallmatrix}\right]$  of a structured polynomial matrix $\left[ \begin{smallmatrix}T(\lambda) & U(\lambda)\\ V(\lambda) & W(\lambda)\end{smallmatrix}\right]$  \\
 \texttt{\bfseries spm2lps}    & Building a pencil based structured linearization $\left[ \begin{smallmatrix}A-\lambda E & B-\lambda F\\ C-\lambda G  & D-\lambda H\end{smallmatrix}\right]$ of a structured polynomial matrix $\left[ \begin{smallmatrix}T(\lambda) & U(\lambda)\\ V(\lambda) & W(\lambda)\end{smallmatrix}\right]$  \\
  \texttt{\bfseries pmkstruct}      & Determination of the Kronecker structure and the multiplicities of infinite poles and zeros  using companion form based linearizations   \\
  \texttt{\bfseries pmeigvals}    & Computation of the finite and infinite eigenvalues using companion form based linearizations  \\
  \texttt{\bfseries pmzeros}    & Computation of the finite and infinite zeros using companion form based linearizations  \\
  \texttt{\bfseries pmzeros1}    & Computation of the finite and infinite zeros using pencil based structured linearizations  \\
  \texttt{\bfseries pmzeros2}    & Computation of the finite and infinite zeros using structured linearizations  \\
  \texttt{\bfseries pmroots}     & Computation of the roots of the determinant of a regular polynomial matrix (i.e., finite zeros)\\
  \texttt{\bfseries pmpoles}     & Computation of the finite and infinite poles using companion form based linearizations \\
  \texttt{\bfseries pmpoles1}  & Computation of the finite and infinite poles using pencil based structured linearizations  \\
  \texttt{\bfseries pmpoles2}  & Computation of the finite and infinite poles using structured linearization\\
  \texttt{\bfseries pmrank}   & Determination of the normal rank of a polynomial matrix\\
  \texttt{\bfseries ispmregular}   & Checking the regularity of a polynomial matrix \\
  \texttt{\bfseries ispmunimodular}    & Checking the unimodularity of a polynomial matrix \\
  \hline
\end{longtable}

\newpage
The following table lists the main functions  for rational matrices available in \textsf{Release v1.1} of the  \textsf{MatrixPencils} package:
%\begin{table}[h]
\begin{longtable}{lp{10.5cm}} Function & Description \\ \hline
  \texttt{\bfseries rm2lspm}     & Construction of a representation of a rational matrix as the sum of its strictly proper part (realized as a structured linearization) and its polynomial part\\
  \texttt{\bfseries rmeval}     & Evaluation of a rational matrix for a given value of its argument\\
  \texttt{\bfseries rm2ls}    & Building a structured linearization $\left[ \begin{smallmatrix}A-\lambda E & B\\ C & D\end{smallmatrix}\right]$ of a rational matrix \\
  \texttt{\bfseries ls2rm}     & Computation of the rational matrix from its structured linearization\\
  \texttt{\bfseries rm2lps}      & Building a pencil based structured linearization
  $\left[ \begin{smallmatrix}A-\lambda E & B-\lambda F\\ C-\lambda G  & D-\lambda H\end{smallmatrix}\right]$
   of a rational matrix\\
  \texttt{\bfseries lps2pm}      & Computation of the rational matrix from its pencil based structured linearization  \\
  \texttt{\bfseries lpmfd2ls}    & Building a structured linearization $\left[ \begin{smallmatrix}A-\lambda E & B\\ C & D\end{smallmatrix}\right]$  of a left polynomial matrix fractional description $D^{-1}(\lambda)N(\lambda)$   \\
  \texttt{\bfseries rpmfd2ls}    & Building a structured linearization $\left[ \begin{smallmatrix}A-\lambda E & B\\ C & D\end{smallmatrix}\right]$  of a right polynomial matrix fractional description $N(\lambda)D^{-1}(\lambda)$   \\
 \texttt{\bfseries lpmfd2lps}    & Building a pencil based structured linearization $\left[ \begin{smallmatrix}A-\lambda E & B-\lambda F\\ C-\lambda G  & D-\lambda H\end{smallmatrix}\right]$ of a left polynomial matrix fractional description $D^{-1}(\lambda)N(\lambda)$   \\
 \texttt{\bfseries rpmfd2lps}    & Building a pencil based structured linearization $\left[ \begin{smallmatrix}A-\lambda E & B-\lambda F\\ C-\lambda G  & D-\lambda H\end{smallmatrix}\right]$ of a right polynomial matrix fractional description $N(\lambda)D^{-1}(\lambda)$   \\
 \texttt{\bfseries rpmfd2lps}    & Building a pencil based structured linearization $\left[ \begin{smallmatrix}A-\lambda E & B-\lambda F\\ C-\lambda G  & D-\lambda H\end{smallmatrix}\right]$ of the inverse of a polynomial matrix   \\
  \texttt{\bfseries pminv2ls}    & Building a structured linearization $\left[ \begin{smallmatrix}A-\lambda E & B\\ C & D\end{smallmatrix}\right]$  of the inverse of a polynomial matrix    \\
  \texttt{\bfseries rmkstruct}      & Determination of the Kronecker structure and the multiplicities of infinite poles and zeros  using companion form based linearizations   \\
  \texttt{\bfseries rmzeros}    & Computation of the finite and infinite zeros using using structured descriptor system linearizations  \\
  \texttt{\bfseries rmzeros1}    & Computation of the finite and infinite zeros using pencil based structured linearizations  \\
  \texttt{\bfseries rmpoles}     & Computation of the finite and infinite poles using using structured descriptor system linearizations \\
  \texttt{\bfseries rmpoles1}  & Computation of the finite and infinite poles using pencil based structured linearizations  \\
  \texttt{\bfseries rmrank}   & Determination of the normal rank of a polynomial matrix\\
\end{longtable}

%\caption[Functions in the \textsf{MatrixPencils}]{Functions in the \textsf{MatrixPencils} package for manipulating polynomial matrices} \label{tab:m}
%\end{table}

\section{Examples}
\label{example}
To illustrate the main concepts related to polynomial and rational matrices, we present two examples which possess all discussed essential structural features and can be handled both analytically and numerically.

\subsection{Example 1} This example, taken from \cite{Door83a},  is a $p\times m$ matrix with $p = m = 3$
\be\label{Ex1} P(\lambda) = \ba{ccc} \lambda^2+\lambda+1 & 4\lambda^2+3\lambda+2 & 2\lambda^2-2 \\
\lambda & 4\lambda -1 & 2\lambda -2 \\
\lambda^2 & 4\lambda^2-\lambda & 2\lambda^2-2\lambda \ea \ee
of degree $d = 2$ and rank $r = 2$. $P(\lambda)$ can be alternatively expressed as the matrix polynomial $P(\lambda) = P_0+P_1\lambda + P_2\lambda^2$ in the standard monomial basis, with
\[ P_0 = \ba{rrr} 1 & 2 & -2 \\ 0 & -1 & -2\\ 0 & 0 & 0 \ea, \qquad P_1 = \ba{rrr} 1 & 3 & 0\\1 & 4 & 2 \\0 & -1 & -2 \ea,
\qquad P_2 = \ba{ccc} 1 & 4 & 2\\ 0 & 0 & 0\\ 1 & 4 & 2 \ea .\]

To study the finite eigenvalue structure we use the unimodular matrices from \cite{Door83a}
\be\label{UV} U(\lambda) = \ba{rcr} 1 & -1 & -1 \\ -\lambda & \lambda+1 & \lambda\\ 0 & -\lambda & 1 \ea, \qquad
V(\lambda) = \ba{rrr} 1 & -3 & 6\\ 0 & 1& -2\\ 0 & 0 & 1\ea , \ee
to obtain the Smith-form of $P(\lambda)$ as
\be\label{SmithEx1} D(\lambda) = U(\lambda)P(\lambda)V(\lambda) = \ba{ccc} 1 & 0 & 0 \\ 0 & \lambda - 1 & 0 \\ 0 & 0 & 0 \ea ,\ee
which exhibits the finite eigenvalue $\lambda = 1$ with partial multiplicities $(0,1)$,  and the (normal) rank $r = 2$ of $P(\lambda)$.

To study the structure at infinity of $P(\lambda)$, we determine the multiplicities of the infinite eigenvalues as the multiplicities of the null eigenvalues of the reversal $P_{rev}(\lambda) := \lambda^2P(1/\lambda)$. We used the following two unimodular matrices
\[ \widetilde U(\lambda) = \ba{crc} -\lambda & 0 & \lambda^2+2\lambda+1 \\ \lambda^2-\lambda+1 & 0 & -\lambda^3-\lambda^2-1\\ 0 & -1 & \lambda \ea, \;\;
\widetilde V(\lambda) = \ba{ccr} 3\lambda+1 & 3\lambda^3+\lambda^2-3\lambda-4 & 6\\ -\lambda & -\lambda^3+\lambda+1& -2\\ 0 & 0 & 1\ea , \]
to obtain the Smith-form of $P_{rev}(\lambda)$ as
\[ D_{rev}(\lambda) = \widetilde U(\lambda)P_{rev}(\lambda)\widetilde V(\lambda) = \ba{ccc} 1 & 0 & 0 \\ 0 & \lambda^2(\lambda - 1) & 0 \\ 0 & 0 & 0 \ea , \]
which exhibits the finite eigenvalue $\lambda = 0$ with partial multiplicities $(0, 2)$, and therefore two infinite eigenvalues with the same partial multiplicities, and, additionally, the finite zero $\lambda = 1$ with partial multiplicities $(0,1)$. This nonzero finite zero of $P_{rev}(\lambda)$ is the reciprocal of the finite zero of $P(\lambda)$, and has evidently the same structural indices. As expected, the (normal) rank  $P_{rev}(\lambda)$ is 2. It follows that the spectrum of $P(\lambda)$, formed of the finite and infinite eigenvalues, is $\mathcal{E} = \{ 1, \infty \}$, with partial multiplicities $(0,1)$ and $(0,2)$, respectively.    The  finite zero structure (according to McMillan \cite{Rose70}) and the finite eigenvalue structure at $\lambda = 1$ coincide, with the finite structural indices  $(0,1)$. The infinite pole-zero structure is given by the infinite structural indices $(-2,0)$ and indicates an infinite pole of multiplicity two and no infinite zero (recall that the partial multiplicities of infinite eigenvalues are in excess with $d = 2$).

If we regard $P(\lambda)$ as a rational matrix, then we can alternatively use for our analysis the local Smith-McMillan form of $P(\lambda)$ (as in \cite{Door83a}). For the finite eigenvalue structure, the analysis based on the Smith form is satisfactory. For the analysis of the infinite structure, we employ two matrices $U_\infty(\lambda)$ and $V_\infty(\lambda)$, which are regular at $\lambda = \infty$,  to
determine the structure of $P(\lambda)$ at infinity. For reference purposes we give the expressions of these matrices
\[ U_\infty(\lambda) =
\ba{ccc} -\frac{1}{\lambda} & 0 & \frac{{\left(\lambda + 1\right)}^2}{\lambda^2}\\ -\frac{\lambda^2 - \lambda + 1}{\lambda\, \left(\lambda - 1\right)} & 0 & \frac{\lambda^3 + \lambda + 1}{\lambda^2\, \left(\lambda - 1\right)}\\ 0 & -1 & \frac{1}{\lambda} \ea,  \qquad
V_\infty(\lambda) =
\ba{ccc} \frac{\lambda+3}{\lambda}  & \frac{ - 4\, \lambda^3 - 3\, \lambda^2 + \lambda + 3}{\lambda^3} & 6\\ -\frac{1}{\lambda} & \frac{\lambda^3 + \lambda^2 - 1}{\lambda^3} & -2\\ 0 & 0 & 1 \ea
\]
and the resulting local Smith-McMillan form at $\infty$
\be\label{SmithInfEx1}  D_\infty(\lambda) = U_\infty(\lambda)P(\lambda)V_\infty(\lambda) = \ba{ccc} (1/\lambda)^{-2} & 0 & 0 \\ 0 & 1 & 0 \\ 0 & 0 & 0 \ea .\ee
The above form shows that $\infty$ is indeed a pole of multiplicity two of $P(\lambda)$ and $P(\lambda)$ has no zeros at infinity.

The last column of $V(\lambda)$ in (\ref{UV}) is a right annihilator of $P(\lambda)$ of degree 0 and, represents a minimal polynomial basis of the right nullspace $\mathcal{N}_r(P)$ of $P(\lambda)$. Similarly, the last row of $U(\lambda)$ in (\ref{UV}) is a left annihilator of $P(\lambda)$ of degree 1 and, represents a minimal polynomial basis of the left nullspace $\mathcal{N}_l(P)$ of $P(\lambda)$. It follows, that the singularity of $P(\lambda)$ is characterized by the right Kronecker index $\epsilon_1 = 1$ and the left Kronecker index $\eta_1 = 0$.

In what follows, we determine the structural properties of $P(\lambda)$ by employing the three types of discussed linearizations.

\subsubsection{Using a companion form linearization}
Using the first Frobenius companion form linearization of $P(\lambda)$ of grade $k = d$, we obtain the pencil $C_1(\lambda) = M_1 - \lambda N_1$, with
\[ M_1 = \ba{rrrrrr}
 -1 & -3 &  0 &  -1 &  -2 &   2 \\
 -1 & -4 & -2 &   0 &   1 &   2 \\
  0 &  1 &  2 &   0 &   0 &   0 \\
  1 &  0 &  0 &   0 &   0 &   0 \\
  0 &  1 &  0 &   0 &   0 &   0 \\
  0 &  0 &  1 &   0 &   0 &   0
  \ea,  \qquad
N_1 = \ba{rrrrrr}
 1 &  4 &  2 &  0 &  0 & 0 \\
 0 &  0 &  0 &  0 &  0 & 0 \\
 1 &  4 &  2 &  0 &  0 & 0 \\
 0 &  0 &  0 &  1 &  0 & 0 \\
 0 &  0 &  0 &  0 &  1 & 0 \\
 0 &  0 &  0 &  0 &  0 & 1
  \ea.
 \]

 The computation of the Kronecker structure of $C_1(\lambda)$ reveals the following: a finite eigenvalue 1 with the partial multiplicities $( 0, 1)$ (not explicitly determined) and hence $\delta_{fin}(C_1) = 1$;  two infinite eigenvalues and the corresponding partial multiplicities $( 0, 2)$, hence $\delta_{\infty}(C_1) = 2$; the right Kronecker index
$\widetilde\epsilon_1 = 1$ and the left Kronecker index $\eta_1 = 1$, and hence $\mu(C_1) = 2$. From this information, we can recover the right Kronecker index $\epsilon_1$ of $P(\lambda)$ (see Proposition \ref{PropCF1}) as $\epsilon_1 = \widetilde\epsilon_1 - (d-1) = 0$. The normal rank of $P(\lambda)$ results as $\rank P(\lambda) = \rank C_1(\lambda)- m(d-1) = 2$, where $\rank C_1(\lambda) = \delta_{fin}(C_1)+\delta_{\infty}(C_1)+\mu(C_1) = 5$.
The finite zero structure of $P(\lambda)$ is the same as the finite eigenvalue structure of $C_1(\lambda)$, while the infinite pole-zero structure of $P(\lambda)$ results from the resulted partial multiplicities of infinite eigenvalues (i.e., (0,2)), which exceed with $d = 2$ the infinite structural indices $(-2,0)$, thus indicating an infinite pole of multiplicity 2 and no infinite zeros.

 Similar results can be obtained
 using the second Frobenius companion form linearization of   $P(\lambda)$.

\subsubsection{Using pencil based linearization}

A strongly minimal realization of $P(\lambda)$ can be determined by inspection, observing that the coefficient matrix $P_2$
of $\lambda^2$ can be expressed in a full rank factorized form  $P_2 = L R$, with
\[   L = [\, 1 \; 0 \; 1 \,]^T, \qquad  R =  [\, 1 \; 4 \; 2 \,] , \]
which immediately leads to the strongly minimal realization of order 1 with
\[ A = -1, \quad E = 0, \quad B = [\, 0 \; 0 \; 0 \,], \quad F = R, \quad C = [\, 0 \; 0 \; 0 \,]^T, \quad G = L, \quad D = P_0, \quad
H = -P_1   .\]

The computation of the Kronecker structure of the system matrix pencil $S(\lambda)$ in (\ref{sysmat}) reveals the following:
a finite eigenvalue 1 with the partial multiplicities $( 0, 1)$  (not explicitly determined) and hence $\delta_{fin}(S) = 1$;  an infinite eigenvalue and the corresponding partial multiplicities $( 0, 1)$, hence $\delta_{\infty}(S) = 1$; the right Kronecker index
$\epsilon_1 = 0$ and the left Kronecker index $\eta_1 = 1$, and hence $\mu(S) = 1$. The system matrix $S(\lambda)$ has a finite zero at 1 and no infinite zeros, and therefore the zero and singular structures of $P(\lambda)$ and $S(\lambda)$ coincide and $\delta^z(P) = 1$.

The computation of the Kronecker structure of the pole pencil $S_p(\lambda)$ in (\ref{TFM-poles}) reveals the following:
no finite eigenvalues and hence $\delta_{fin}(S_p) = 0$;  seven infinite eigenvalues and the corresponding partial multiplicities $( 1, 1, 1, 1, 3)$, hence $\delta_{\infty}(S_p) = 7$; no right and left Kronecker indices, and hence $\mu(S_p) = 0$. The pole pencil $S_p(\lambda)$ has no finite zeros and has two infinite zeros of multiplicity two, and therefore the pole and singular structures of $P(\lambda)$ and the zeros and singular structures of $S_p(\lambda)$ coincide and $\delta^p(P) = \delta^z(S_p) = 2$.
The condition $\delta^p(P) = \delta^z(P)+\mu(P)$ is fulfilled, because  $\delta^p(P) = \delta^z(S_p) = \delta^z(S) + \mu(S)$.

From the knowledge of the infinite pole-zero structure with the infinite structural indices $(-2,0)$ (i.e., two infinite poles and no infinite zero), we can determine the infinite eigenvalue structure by shifting these values with $d = 2$ (the degree of $P(\lambda)$). We obtain the expected partial multiplicities of infinite eigenvalues $(0,2)$.

\subsubsection{Using descriptor system based linearizations}

A third possibility to determine the  pole- zero and the singular Kronecker structures is to use a descriptor system realization based linearization of $P(\lambda)$ of the form (\ref{sysmatD})
where $(A-\lambda E, B, C, D)$ satisfies (\ref{TFM2}).
Recall that, if the descriptor realization is irreducible (i.e., controllable and observable), then the
zero and singular structures of $S(\lambda)$ in (\ref{sysmatD}) and $P(\lambda)$ coincide.

Consider the irreducible realization of order $n = 4$ (also used in \cite{Misr94}) with
\[
A = \ba{rrrr}
0 & 0 &  1 &  0 \\
0 & 0 &  0 &  1 \\
1 & 0 &  0 &  0 \\
0 & 1 & -1 &  0\ea, \quad
E = \ba{rrrr}
1 &  0 &  0 &  0 \\
0 &  1 &  0 &  0 \\
0 &  0 &  0 &  0 \\
0 &  0 &  0 &  0\ea, \quad
B = \ba{rrr}
0 &  0 &  0 \\
0 &  0 &  0 \\
1 &  4 &  2 \\
0 & -1 & -2\ea, \]
\[ C = \ba{rrrr}
0 &  0 &  -1 &  -1  \\
0 &  0 &  -1 &  0  \\
0 &  0 &  0  & -1\ea, \quad
D = \ba{rrr}
1 &  2  & -2  \\
0 &  -1 &  -2 \\
0 &  0  & 0\ea .\]

The computation of the Kronecker structure of the system matrix pencil $S(\lambda)$ in (\ref{sysmatD}) reveals the following:
a finite eigenvalue 1 with the partial multiplicities $( 0, 1)$ (not explicitly determined) and hence $\delta_{fin}(S) = 1$;  four infinite eigenvalues and the corresponding partial multiplicities $( 1, 1, 1, 1)$, hence $\delta_{\infty}(S) = 4$; the right Kronecker index
$\epsilon_1 = 0$ and the left Kronecker index $\eta_1 = 1$, and hence $\mu(S) = 1$. The system matrix $S(\lambda)$ has a finite zero at 1 and no infinite zeros, and therefore the zero and singular structures of $P(\lambda)$ and $S(\lambda)$ coincide.
%From (\ref{stronginf}), we obtain $\delta_\infty(P) = 2$, thus the eigenvalues of $P(\lambda)$ are $\{ 1, \infty, \infty \}$.

The computation of the Kronecker structure of the pole pencil $S_p(\lambda) = A-\lambda E$ in (\ref{TFM-poles}) reveals the following:
no finite eigenvalues and hence $\delta_{fin}(S_p) = 0$;  four infinite eigenvalues and the corresponding partial multiplicities $( 1, 3)$, hence $\delta_{\infty}(S_p) = 4$; no right and left Kronecker indices, and hence $\mu(S_p) = 0$. The pole pencil $S_p(\lambda)$ has no finite zeros and has two infinite zeros of multiplicity two, and therefore the pole and singular structures of $P(\lambda)$ and the zeros and singular structures of $S_p(\lambda)$ coincide and the condition $\delta^p(P) = \delta^z(P)+\mu(P)$ is fulfilled, because  $\delta^p(P) = \delta^z(S_p) = \delta^z(S) + \mu(S) = 2$.

Similar results have been obtained using a minimal descriptor realization (i.e., without non-dynamic modes) of order $n = 3$. For reference purposes we give the matrices of employed realization
\[
A = \ba{rrr}
0 & 1 &  0 \\
0 &  0 &  1 \\
1 & 0 &  0 \ea, \quad
E = \ba{rrr}
1 &  0 &  0  \\
0 &  1 &  0  \\
0 &  0 &  0 \ea, \quad
B = \ba{rrr}
0 & -1 & -2 \\
0 &  0 &  0 \\
1 &  4 &  2 \ea, \]
\[ C = \ba{rrrr}
0 &  -1 &  -1  \\
0 &  -1 &  0  \\
0 &  0  & -1\ea, \quad
D = \ba{rrr}
1 &  3  & 0  \\
0 &  0 &  0 \\
0 &  0  & 0\ea . \]
It is worth mentioning, that for the computation of zeros and poles, the use of a minimal realization instead of an irreducible one has no practical advantages. This is because the determination of a minimal realization usually involves, besides the determination of an irreducible realization using orthogonal similarity transformations, the additional step of eliminating the non-dynamic modes, which however involves matrix inversions and thus cannot be performed using only orthogonal reductions.

 \subsection{Example 2} This example is a rational matrix $R(\lambda)$ defined as $R(\lambda) = \frac{1}{\lambda+1}P(\lambda)$, where $P(\lambda)$ is the polynomial matrix (\ref{Ex1}) employed in Example 1. For completeness, we give below the explicit expression of $R(\lambda)$
 \[ R(\lambda) = \left[\begin{array}{ccc} \displaystyle\frac{\lambda^2 + \lambda + 1}{\lambda + 1} & \displaystyle\frac{4\, \lambda^2 + 3\, \lambda + 2}{\lambda + 1} & 2\, \lambda - 2\\ \\[-3mm] \displaystyle\frac{\lambda}{\lambda + 1} & \displaystyle\frac{4\, \lambda - 1}{\lambda + 1} & \displaystyle\frac{2\, \lambda - 2}{\lambda + 1}\\ \\[-3mm] \displaystyle\frac{\lambda^2}{\lambda + 1} & \displaystyle\frac{4\, \lambda^2 - \lambda}{\lambda + 1} & \displaystyle\frac{2\, \lambda^2 -  2\, \lambda}{\lambda + 1} \end{array}\right].
\]

 For the analysis of the finite pole-zero structure, we can alternatively determine the (non-local) Smith-McMillan form \cite{Rose70} directly from the Smith-form (\ref{SmithEx1}) of $P(\lambda)$ as
 \[ \frac{1}{\lambda+1} D(\lambda) = \left[\begin{array}{ccc} \displaystyle\frac{1}{\lambda + 1} & 0 & 0\\ 0 & \displaystyle\frac{\lambda - 1}{\lambda + 1} & 0\\ 0 & 0 & 0 \end{array}\right]
\]
from which we can read out the presence of two simple poles at $\lambda = -1$ and of a zero at $\lambda = 1$. For the analysis of the infinite pole-zero structure, we can adjust the local Smith-McMillan form at infinity in (\ref{SmithInfEx1}), by computing
\[ \frac{\lambda}{\lambda+1}U_\infty(\lambda)R(\lambda)V_\infty(\lambda) = \ba{ccc} (1/\lambda)^{-1} & 0 & 0 \\ 0 & 1/\lambda & 0 \\ 0 & 0 & 0 \ea .\]
This shows that $\lambda = \infty$ is both a pole as well as a zero of $R(\lambda)$.

The last column of $V(\lambda)$ in (\ref{UV}) is a right annihilator of $R(\lambda)$ of degree 0 and, represents a minimal polynomial basis of the right nullspace $\mathcal{N}_r(R)$ of $R(\lambda)$. Similarly, the last row of $U(\lambda)$ in (\ref{UV}) is a left annihilator of $R(\lambda)$ of degree 1 and, represents a minimal polynomial basis of the left nullspace $\mathcal{N}_l(R)$ of $R(\lambda)$. It follows, that the singularity of $R(\lambda)$ is characterized by the right Kronecker index $\epsilon_1 = 1$ and the left Kronecker index $\eta_1 = 0$.

In what follows, we determine the structural properties of $R(\lambda)$ by employing system matrix pencils based linearizations.

\subsubsection{Using pencil based linearization}

A strongly minimal realization of $R(\lambda)$ can be determined from the additive decomposition
\be\label{Rsppol} R(\lambda) = R_{sp}(\lambda) + R_{pol}(\lambda) \ee
with the strictly proper part
\[  R_{sp}(\lambda) = \left[\begin{array}{ccc} \displaystyle\frac{1}{\lambda + 1} & \displaystyle\frac{3}{\lambda + 1} & 0\\ \\[-3mm] \displaystyle\frac{-1}{\lambda + 1} & \displaystyle\frac{-5}{\lambda + 1} & \displaystyle\frac{-4}{\lambda + 1}\\ \\[-3mm] \displaystyle\frac{1}{\lambda + 1} & \displaystyle\frac{5}{\lambda + 1} & \displaystyle\frac{4}{\lambda + 1} \end{array}\right]\]
and polynomial part
\[ R_{pol}(\lambda) = \ba{ccc} \lambda & 4\lambda-1 & 2\lambda -2 \\ 1 & 4 & 2 \\ \lambda -1 & 4\lambda -5 & 2\lambda-4 \ea .\]
For the strictly proper part a non-minimal realization of order 8 can be easily constructed by inspection in the form (\ref{Rspnonmin}) using individual realizations of entries. Also, controllable or observable realizations of order 3 in the form (\ref{spcontr}) can be easily constructed by inspection. From any of these realizations, a second order minimal realization can be obtained by removing the uncontrollable and/or unobservable eigenvalues using orthogonal transformations based algorithms as in \cite{Door81}. Alternatively, a minimal realization can be obtained using a full rank factorization of $R_{sp}(\lambda)$, which can be constructed in the form
\be\label{Rsp} R_{sp}(\lambda) = L \widetilde R_{sp}(\lambda) , \ee
where
\[ L = \ba{rr} 1 & 1\\ 0 & 1\\ 0 & -1 \ea, \quad \widetilde R_{sp}(\lambda) = \left[\begin{array}{ccc} \displaystyle\frac{1}{\lambda + 1} & \displaystyle\frac{3}{\lambda + 1} & 0\\ \\[-3mm] \displaystyle\frac{-1}{\lambda + 1} & \displaystyle\frac{-5}{\lambda + 1} & \displaystyle\frac{-4}{\lambda + 1} \end{array}\right] .\]
A second order observable realization of $\widetilde R_{sp}(\lambda)$ is obtained by inspection in the dual form of  (\ref{spcontr}). This realization is also controllable, thus it is minimal. Overall, we obtain for $R_{sp}(\lambda)$ in (\ref{Rsp}) a minimal realization
$(A_{sp}-\lambda I,B_{sp},C_{sp},0)$ with the matrices
\be\label{minsp1} A_{sp} = \ba{rr} -1 & 0 \\ 0 & -1 \ea, \quad
B_{sp} = \ba{rrr}  1 & 3 &  0\\
 -1 & -5 & -4 \ea, \quad
 C_{sp} = \ba{rr} 1 & 1\\ 0 & 1\\ 0 & -1 \ea  .\ee

For the polynomial part a zeroth-order pencil based realization is given by
\be\label{Rpol1} R_{pol}(\lambda) = D_{pol}-\lambda H_{pol}, \ee
with
\[ D_{pol} = \ba{rrr} 0 & -1 & -2\\ 1 & 4 & 4\\ -1 & -5 &  -4 \ea, \quad H_{pol} = \ba{rrr} -1 & -4 & -2 \\ 0 & 0 & 0 \\ -1 & -4 & -2 \ea .\]
The overall pencil-based realization is $(A_{sp}-\lambda I,B_{sp},C_{sp},D_{pol}-\lambda H_{pol})$.
The zeros and Kronecker structure of $R(\lambda)$ are those of the system matrix
\[ \overline S(\lambda) = \ba{cc} A_{sp}-\lambda I & B_{sp} \\
C_{sp} & D_{pol}- \lambda H_{pol} \ea \, ,\]
while the finite poles are the eigenvalues of $A_{sp}$ and the infinite poles are the infinite zeros of the subpencil in (\ref{Shat})
\[ \widehat S_\infty(\lambda) := \ba{cc}   0 & I_p \\ I_m & 0\ea - \lambda  \ba{ccc}  H_{pol} & 0 \\ 0 & 0  \ea .\]

The computation of the zeros and of Kronecker structure of the system matrix pencil $\overline S(\lambda)$  revealed the following: a finite zero at 1 with the partial multiplicities $( 0, 1)$  (not explicitly determined) and an infinite eigenvalue with partial multiplicities $( 0, 2)$, which corresponds to an infinite zero; the right Kronecker index
$\epsilon_1 = 0$ and the left Kronecker index $\eta_1 = 1$, and hence $\mu(\overline S) = 1$. The system matrix $ \overline S(\lambda)$ has a finite zero at 1 and an infinite zero, and therefore the zero and singular structures of $R(\lambda)$ and $ \overline S(\lambda)$ coincide and $\delta^z(R) = 2$ and $\mu(R) = 1$.

The computation of the poles as the zeros of the pole pencil (\ref{TFM-poles}) has been performed by exploiting its particular structure, such that the finite poles are formed of two poles ar $-1$ (these are the eigenvalues of $A_{sp}$) and an infinite pole (which is an infinite zero of $\widehat S_\infty(\lambda)$ above). Since $\delta^p(R) = 3$, we have the condition (\ref{IST-VD}) fulfilled.

\subsubsection{Using descriptor system based linearization}
An irreducible descriptor system realization of $R(\lambda)$ can be also determined from the additive decomposition (\ref{Rsppol}). For the strictly proper part we already determined a minimal realization $(A_{sp}-\lambda I,B_{sp},C_{sp},0)$ of order 2 with the matrices given in (\ref{minsp1}). For the polynomial part $R_{pol}(\lambda)$ a controllable or observable descriptor system realizations of the form  $(I-\lambda E_{pol}, B_{pol}, C_{pol}, D_{pol})$ can be constructed as in (\ref{sscon}) or (\ref{ssobs}), respectively. Both realizations have order 6 and irreducible realizations of order 2 can be obtained by removing the unobservable or uncontrollable infinite eigenvalues using orthogonal transformation based methods proposed in \cite{Door81}.
The simple form of $R_{pol}(\lambda)$ in (\ref{Rpol1}) allows to directly obtain a second order realization of the term $\lambda H_{pol}$, observing that
\[  \lambda H_{pol} = \ba{r} 1\\ 0 \\ 1 \ea \lambda \; [\, -\!1 \; -\!4 \; -\!2\,] \]
Using for $\lambda$ the irreducible realization $\left(\left[\begin{smallmatrix} 1 & -\lambda \\ 0 & 1 \end{smallmatrix}\right], \left[\begin{smallmatrix} 0 \\ -1 \end{smallmatrix}\right], [\, 1 \; 0\,], 0 \right)$, we obtain the matrices of the realization of $R_{pol}(\lambda)$ as
\[ {\arraycolsep = 0.7mm A_{pol}-\lambda E_{pol} = \ba{rr} 1 & -\lambda \\ 0 & 1 \ea, \quad B_{pol} = \ba{rrr} 0 & 0 & 0 \\-1 & -4 & -2 \ea, \quad C_{pol} = \ba{rr} 1 & 0\\0 & 0\\ 1 & 0 \ea, \quad D_{pol} = \ba{rrr} 0 & -1 & -2\\ 1 & 4 & 4\\ -1 & -5 &  -4 \ea }. \]

The overall descriptor system realization has the form (\ref{Rtot}). The zeros and Kronecker structure of $R(\lambda)$ are those of the system matrix
\[ \overline S(\lambda) =  \ba{ccc} A_{sp}-\lambda I & 0 & B_{sp} \\ 0 & I-\lambda E_{pol} & B_{pol} \\
C_{sp} & C_{pol} & D_{pol} \ea \, ,\]
while the poles are the zeros (finite and infinite) of the pole pencil 
\[ \widehat S(\lambda) = \ba{cc} A_{sp}-\lambda I & 0 \\ 0 & I-\lambda E_{pol} \ea .\] The  finite poles are therefore the eigenvalues of $A_{sp}$ and the infinite poles are the infinite zeros of the subpencil $I - \lambda E_{pol}$.

The computed pole-zero structure of $R(\lambda)$ is identical to that computed with the pencil-based linearization approach.

\section{Conclusions}
In this article we presented the main theoretical results which are relevant for the determination of the Kronecker and pole-zero structures of polynomial and rational matrices using linearization based computational techniques. The companion form based linearizations served as basis to implement the basic structural analysis functions of the \textsf{MatrixPencils} package to compute eigenvalues, singular structures and pole-zero structures of polynomial matrices. Alternatively, linearizations based on pencil and descriptor system representations are used to implement functions for the determination of the pole-zero and singular structures of both polynomial and rational matrices.

Some useful links for the \textsf{MatrixPencils} package are listed below:
\begin{itemize}
\item[--] download site of the latest release \url{https://github.com/andreasvarga/MatrixPencils.jl};
\item[--] alternative download site \url{https://zenodo.org/record/4004252};
\item[--] latest version of the documentation \url{https://andreasvarga.github.io/MatrixPencils.jl/dev/};
\item[--] complete list of available functions \url{https://sites.google.com/site/andreasvargacontact/home/software/matrix-pencils-in-julia}.
\end{itemize}

%\newpage
\bibliographystyle{siamplain}
%\bibliography{../../varga}

\end{document}